\documentclass{amsart}
\usepackage{tabu}
\usepackage[margin=95pt]{geometry}
\usepackage{enumitem}
\allowdisplaybreaks

\usepackage{amssymb}
\usepackage{amsmath}
\usepackage{amscd}
\usepackage{amsbsy}
\usepackage{comment}
\usepackage[matrix,arrow]{xy}
\usepackage{hyperref}
\DeclareMathOperator{\Rad}{Rad}
\DeclareMathOperator{\LCM}{LCM}

\DeclareMathOperator{\Norm}{Norm}
\DeclareMathOperator{\Disc}{Disc}

\DeclareMathOperator{\Gal}{Gal}

\DeclareMathOperator{\ord}{ord}

\newcommand{\Q}{{\mathbb Q}}
\newcommand{\Z}{{\mathbb Z}}

\newcommand{\F}{{\mathbb F}}
\newcommand{\cA}{\mathcal{A}}
\newcommand{\cB}{\mathcal{B}}

\newcommand{\cL}{\mathcal{L}}

\newcommand{\cE}{\mathcal{E}}
\newcommand{\cF}{\mathcal{F}}
\newcommand{\cG}{\mathcal{G}}
\newcommand{\cH}{\mathcal{H}}

\newcommand{\cC}{\mathcal C}
\newcommand{\cD}{\mathcal D}

\newcommand{\OO}{{\mathcal O}}

\newcommand{\fp}{\mathfrak{p}}
\newcommand{\fq}{\mathfrak{q}}
\newcommand{\fa}{\mathfrak{a}}
\newcommand{\fb}{\mathfrak{b}}

\def\mod#1{{\ifmmode\text{\rm\ (mod~$#1$)}
\else\discretionary{}{}{\hbox{ }}\rm(mod~$#1$)\fi}}

\begin {document}

\newtheorem{thm}{Theorem}
\newtheorem{lem}{Lemma}[section]
\newtheorem{prop}[lem]{Proposition}

\newtheorem{cor}[lem]{Corollary}

\theoremstyle{definition}

\theoremstyle{remark}

\title[Ramanujan tau function]
{Odd Values of the Ramanujan tau function}
\author[Michael Bennett]{Michael A. Bennett}
\address{Department of Mathematics, University of British Columbia, Vancouver, B.C., V6T 1Z2 Canada}
\email{bennett@math.ubc.ca}

\author{Adela Gherga}
\address{Mathematics Institute, University of Warwick, Coventry CV4 7AL, United Kingdom}
\email{Adela.Gherga@warwick.ac.uk}

\author{Vandita Patel}
\address{Department of Mathematics, University of  Manchester, Manchester M13 9PL, United Kingdom}
\email{vandita.patel@manchester.ac.uk}

\author{Samir Siksek}
\address{Mathematics Institute, University of Warwick, Coventry CV4 7AL, United Kingdom}
\email{S.Siksek@warwick.ac.uk}
\thanks{ The first-named author is supported by NSERC.
The second and fourth-named authors are
supported by an EPSRC Grant EP/S031537/1 \lq\lq Moduli of
elliptic curves and classical Diophantine problems\rq\rq.}

\date{\today}


\begin {abstract}
We prove a number of results regarding odd values of the 
Ramanujan $\tau$-function.
For example, we prove the existence of an effectively computable positive
constant $\kappa$ such that if $\tau(n)$ is odd and $n \ge 25$
then either
\[
	P(\tau(n)) \; > \; \kappa \cdot \frac{\log\log\log{n}}{\log\log\log\log{n}}
\]
or there exists a prime $p \mid n$ with $\tau(p)=0$. Here $P(m)$ denotes
the largest prime factor of $m$.
We also solve the equation $\tau(n)=\pm 3^{b_1} 5^{b_2} 7^{b_3} 11^{b_4}$
and the equations $\tau(n)=\pm q^b$ where $3\le q < 100$ is prime and the exponents are arbitrary nonnegative integers.
We make use of a variety of methods, including the Primitive Divisor
Theorem of Bilu, Hanrot and Voutier, bounds for solutions
to Thue--Mahler equations due to Bugeaud and Gy\H{o}ry,
and the modular approach via Galois representations of Frey-Hellegouarch
elliptic curves.
\end {abstract}
\maketitle

\section{Introduction}

The {\it Ramanujan $\tau$-function} $\tau (n)$ is defined via the expansion
\begin{equation} \label{taudef}
\Delta (z) = q \prod_{n=1}^\infty (1-q^n)^{24} = \sum_{n=1}^\infty \tau (n) q^n, \; \; \; \; q=e^{2\pi i z}.
\end{equation}
It was conjectured by Ramanujan \cite{Ram} and proved by Mordell \cite{Mor} that $\tau (n)$ is a multiplicative function, i.e. that
$$
\tau (n_1n_2)=\tau (n_1) \tau (n_2),
$$
for all coprime pairs of positive  integers $n_1$ and $n_2$. Further, we have
$$
 \sum_{n=1}^\infty \tau (n) q^n \equiv q \prod_{n=1}^\infty (1+q^{8n})^{3} \equiv q \prod_{n=1}^\infty   (1-q^{8n}) (1+q^{8n})^{2}   \equiv \sum_{n=0}^\infty q^{(2n+1)^2} \mod{2},
 $$
 via Jacobi's triple product formula,
whence $\tau(n)$ is odd precisely when $n$ is an odd square and, in particular, $\tau (p)$ is even for every prime $p$.

 Amongst the many open questions about the possible values of $\tau (n)$, the most notorious is a conjecture of Lehmer \cite{Leh} to the effect that $\tau (n)$ never vanishes.
  In terms of the size of values of $\tau$, one has  the upper bound of Deligne \cite{Deligne} (originally conjectured by Ramanujan) :
 \begin{equation} \label{Del:bound}
 \left| \tau (p) \right| \leq 2 \cdot p^{11/2},
 \end{equation}
 valid for prime $p$.
In the other direction, Atkin and Serre \cite{Ser}  conjectured (as a strengthening of Lehmer's conjecture) that, for $\epsilon > 0$,
$$
 \left| \tau (p) \right| \gg_\epsilon p^{9/2-\epsilon},
 $$
so that, in particular, given a fixed integer $a$, there are at most finitely many primes $p$ for which $\tau (p)=a$.
While this problem remains open, in the special case where the integer $a$ is odd, Murty, Murty and Shorey \cite{MuMuSh} proved that the equation
\begin{equation} \label{equal}
\tau (n)=a,
\end{equation}
has at most finitely many solutions in integers $n$ (note that, in this case, $n$ is necessarily an odd square). More precisely, they demonstrated the existence of an effectively computable positive constant  $c$ such that if $\tau (n)$ is odd, then
$$
|\tau (n)| > (\log (n))^c.
$$
A number of recent papers have treated the problem of explicitly demonstrating that equation (\ref{equal}) has, in fact, no solutions, for various odd values of $a$, including $a=\pm 1$ (Lygeros and Rozier \cite{LyRo}), $|a| < 100$ an odd prime (Balakrishnan, Craig and Ono \cite{BaCrOn}, Balakrishnan, Craig, Ono and Tsai \cite{BaCrOnTs},  Dembner and Jain \cite{DeJa}), and $|a| < 100$ an odd integer (Hanada and Madhukara \cite{HaMa}).

In this paper, we derive what might be considered a non-Archimedean analogue of the work of Murty, Murty and Shorey.
Let us define $P(m)$ to be the greatest prime factor of an integer $|m|>1$. We prove the following.

\begin{thm} \label{thm1}
There exists an effectively computable constant $\kappa > 0$ such that if $\tau (n)$ is odd, with $n \geq 25$, then either
\begin{equation} \label{guppy}
P(\tau (n)) \; > \; \kappa  \cdot \frac{ \log \log  \log n}{\log \log \log \log n},
\end{equation}
or there exists a prime $p \mid n$ for which $\tau (p)=0$.
\end{thm}

Recall that a {\it powerful number} (also known as {\it squarefull} or {\it $2$-full}) is defined to be an integer $n$ with the property that if a prime $p \mid n$, then necessarily $p^2 \mid n$. Equivalently, we can write such an integer as $n=a^2b^3$, where $a$ and $b$ are integers. Our techniques actually show the following (from which Theorem \ref{thm1} is an immediate consequence).
\begin{thm} \label{thm1a}
We have
$$
\lim_{n \rightarrow \infty} P(\tau (n)) = \infty,
$$
where the limit is taken over powerful numbers $n$ for which $\tau (p) \neq 0$ for each $p \mid n$.
More precisely, there exists an effectively computable constant $\kappa > 0$ such that if $n \geq 25$ is powerful,  either
\begin{equation}\label{eqn:lb}
	P(\tau (n)) \; > \; \kappa \cdot \frac{ \log \log  \log n}{\log \log \log \log n}.
\end{equation}
or there exists a prime $p \mid n$ for which $\tau (p)=0$.
\end{thm}

The restriction that $n$ has no prime divisors $p$ for which $\tau (p)=0$ is, in fact, necessary if one wishes to obtain a lower bound upon  $P(\tau (n))$ that tends to $\infty$ with $n$.
Indeed, one may observe that, if $\tau (p) =0$, then
(see \eqref{eqn:recursion} below)
$$
P \left( \tau (p^{2k}) \right) = P \left( (-1)^k p^{11k} \right) =p
$$
is bounded independently of $k$. While Lehmer's conjecture remains unproven, we do know that if there is a prime $p$ for which $\tau (p)=0$, then
\begin{equation}\label{eqn:Derickx}
p > 816212624008487344127999,
\end{equation}
by work of Derickx, van Hoeij and Zeng \cite{DeHoZe}.

Theorem~\ref{thm1a} is an easy consequence of the following
result.
\begin{thm}\label{thm:main}
There is a computable positive constant $\eta$ such that
for any prime $p$ with $\tau(p) \ne 0$ and any $m \ge 2$, 
\begin{equation}\label{eqn:required}
	P(\tau(p^{m})) \; > \; \eta \cdot \frac{ \log \log  (p^{m})}{\log \log \log (p^{m})}.
\end{equation}
\end{thm}
We note that Theorem~\ref{thm:main} implies Theorem~\ref{thm1a}. Indeed, let $n$
be a powerful number and $p^{m}$ be the largest prime power divisor of
$n$.  Then $m \ge 2$, and $p^m \gg \log{n}$, whence \eqref{eqn:lb}  follows  immediately from \eqref{eqn:precise}. Our arguments show the following.
\begin{thm}\label{thm:main2}
Let $m \ge 2$.  There is a  computable positive constant $\delta(m)$, depending
only on $m$, such that
for any prime $p$ with $\tau(p) \ne 0$, 
\begin{equation}\label{eqn:precise}
	P(\tau(p^{m})) \; > \; \delta(m) \cdot \log\log(p). 
\end{equation}
\end{thm}

We note that our bounds neither imply nor are implied by work of Luca and Shparlinski \cite{LuSh} who proved that
$$
P \left( \tau (p) \tau(p^2) \tau(p^3) \right) \gg \frac{\log \log  (p) \log \log \log (p)}{\log \log \log \log (p)}.
$$

\medskip

To demonstrate that these results and the techniques underlying them are somewhat practical, we prove the
 following computational ``coda'', solving equation (\ref{equal}) where the prime divisors of $a$, rather than $a$ itself, are fixed.

\begin{thm} \label{thm2a}
If $n$ is a powerful positive integer, then either $n=8$, where we have
$$
\tau (8) = 2^9 \cdot 3 \cdot 5 \cdot 11,
$$
or
$$
P(\tau (n)) \geq 13.
$$
\end{thm}

\begin{cor} \label{thm2}
If $n$ is a positive integer for which $\tau (n)$ is odd, then
\begin{equation} \label{boundary}
P(\tau (n)) \geq 13.
\end{equation}
In other words, the equation
\begin{equation} \label{fred}
\tau (n) = \pm 3^\alpha 5^\beta 7^\gamma 11^\delta
\end{equation}
has no solutions in integers $n \geq 2$ and $\alpha, \beta, \gamma, \delta \geq 0$.
\end{cor}

It is conjectured that $|\tau (n)|$ takes on infinitely many prime values, the smallest of which corresponds to
$$
\tau (251^2) = -80561663527802406257321747.
$$
Our arguments enable us to eliminate the possibility of powers of small primes arising as values of $\tau$. By way of example, we have the following.
\begin{thm} \label{thm3}
The equation
$$
\tau (n) = \pm q^\alpha
$$
has no solutions in prime $q$ with $3 \leq q < 100$, and $\alpha \ge 0$, $n \ge 2$ integers.
\end{thm}

It is worth observing that the techniques we employ here are readily extended to treat more generally coefficients $\lambda_f(n)$ of cuspidal newforms of (even) weight $k \geq 4$ for $\Gamma_0(N)$, with trivial character and $\lambda_f(p)$ even for suitably large prime $p$; our results correspond to the case of $\Delta (z)$ in (\ref{taudef}), where $k=12$ and $N=1$. For simplicity, we will restrict our attention to $\tau (n)$ and $\Delta (z)$; readers interested in the more general situation should consult the paper of Murty and Murty \cite{MuMu} (see also \cite{BaCrOnTs}).

We should also comment on the particular choice of the constant $``13"$ on the right hand side of inequality (\ref{boundary}) in Corollary \ref{thm2} (and analogously in Theorem \ref{thm2a}). As we shall observe, the weaker result with $13$ replaced by $11$ (corresponding to equation (\ref{fred}) with $\delta=0$) reduces via local arguments to the resolution of a single Thue equation; this is the content of Proposition 6 of  Luca, Maboso, Stanica \cite{LuMaSt}. Corollary \ref{thm2}  as stated requires (apparently at least) the full use of our various techniques, including the Primitive Divisor
Theorem,  solution of a variety of Thue--Mahler equations,
and resolution of hyperelliptic equations through appeal to the modularity of Galois representations attached to Frey-Hellegouarch
elliptic curves. A stronger version of Corollary \ref{thm2}  with $13$ replaced by $17$ in  (\ref{boundary}) is possibly within range of this approach, though computationally significantly more involved. An analogous result with $13$ replaced by $19$  would likely require fundamentally new ideas.

\bigskip

This paper is organized as follows. 
In Section~\ref{sec:congruences}, we recall some standard congruences
for the Ramanujan-tau function that we  use. later in the paper.
In Section~\ref{sec:lucas}, we connect the sequence
$m \mapsto \tau(p^{m-1})$, for a fixed prime $p$, to a Lucas sequence $\{u_m\}$, allowing
us to appeal to  the Primitive Divisor Theorem
of Bilu, Hanrot and Voutier. 
In Section~\ref{sec:Fn}, we introduce a sequence of homogenous
polynomials $\Psi_m(X,Y) \in \Z[X,Y]$ that are intimately connected to the $\{u_m\}$.
We will use these polynomials in Section~\ref{sec:Bugeaud}, together  with  a theorem of Bugeaud on prime divisors of $ax^u+by^v$, to prove Theorem~\ref{thm:main2}.
In Section~\ref{sec:thm1a}, we relate the 
equation $\tau(p^m)=\pm p_1^{\alpha_1} \cdots p_r^{\alpha_1}$
to a Thue--Mahler equation, whence a theorem of Bugeaud and Gy\H{o}ry
enables us to deduce Theorem~\ref{thm:main}.
In Sections~\ref{sec:taup2}, \ref{sec:taup4} and \ref{sec:taup3},
we solve the equations $\tau(p^k)= \pm q^b$ where $k \in \{2, 4 \}$,
$p$ and $q$ are prime, and $3 \le b <100$, and also the equations
$\tau(p^k)= \pm 3^{b_1} 5^{b_2} 7^{b_3} 11^{b_4}$, where $2 \le k \le 4$
and $p$ is prime. Our method in Section~\ref{sec:taup2}, \ref{sec:taup4} and \ref{sec:taup3}
is to associate to a hypothetical solution a Frey-Hellegouarch curve and relate 
this to a weight $2$ modular form
of small level, using recipes of the first author and Skinner
which in turn builds on the modularity of elliptic curves
due to Wiles, Breuil, Conrad Diamond and Taylor, and on
Ribet's Level-Lowering Theorem. In Section~\ref{sec:thm2a},
we prove Theorem~\ref{thm2a}, by combining the results of
Sections~\ref{sec:taup2}, \ref{sec:taup4} and \ref{sec:taup3},
and using the Primitive Divisor Theorem.
Finally, in Section~\ref{sec:thm3}, we prove Theorem~\ref{thm3};
the results of previous sections allow us to
reduce the equation $\tau(n)=\pm q^a$ with $3\le q < 100$
prime to Thue--Mahler equations of high degree, which are then
solved using an algorithm of the second and fourth author, von K\"{a}nel and Matschke.

\section{Congruences for the $\tau$ function} \label{sec:congruences}

For future use, it will be of value for us to record some basic arithmetic facts about $\tau (n)$;
these are taken from Swinnerton--Dyer's article \cite{SwD}.
Here $\sigma_v(n)$ denotes the sum of the $v$-th powers of the divisors
of $n$. 

\begin{equation}\label{mod2}
\begin{cases}
\tau(n) \equiv \sigma_{11}(n) \mod{2^{11}} & \text{if $n \equiv 1 \mod{8}$}\\
\tau(n) \equiv 1217 \cdot \sigma_{11}(n) \mod{2^{13}} & \text{if $n \equiv 3 \mod{8}$}\\
\tau(n) \equiv 1537 \cdot \sigma_{11}(n) \mod{2^{12}} & \text{if $n \equiv 5 \mod{8}$}\\
\tau(n) \equiv 705 \cdot \sigma_{11}(n) \mod{2^{14}} & \text{if $n \equiv 7 \mod{8}$}\\
\end{cases}
\end{equation}

\begin{equation}\label{mod3}
\tau(n) \equiv n^{-610} \cdot \sigma_{1231}(n)
\begin{cases}
\mod{3^6} & \text{if $n \equiv 1 \mod{3}$}\\
\mod{3^7} & \text{if $n \equiv 2 \mod{3}$}\\
\end{cases}
\end{equation}

\begin{equation}\label{mod5}
\tau(n) \equiv n^{-30} \sigma_{71}(n) \mod{5^3} \qquad
\text{if $5 \nmid n$}
\end{equation}

\begin{equation}\label{mod7}
\tau(n) \equiv n \cdot \sigma_9(n)
\begin{cases}
\mod{7} & \text{if $n \equiv 0$, $1$, $2$ or $4 \mod{7}$}\\
\mod{7^2} & \text{if $n \equiv 3$, $5$ or $6 \mod{7}$}\\
\end{cases}
\end{equation}

\begin{equation}\label{mod23}
\begin{cases}
\tau(p) \equiv 0 \mod{23} & \text{if $p$ is a quadratic non-residue
mod $23$}\\
\tau(p) \equiv 2 \mod{23} & \text{if $p=u^2+23v^2$ with $u \ne 0$}\\
\tau(p) \equiv -1 \mod{23} & \text{for other $p \ne 23$}
\end{cases}
\end{equation}

\begin{equation}\label{mod691}
\tau(n) \equiv \sigma_{11}(n) \mod{691}.
\end{equation}

\begin{lem}\label{lem:mod7}
Let $p \ne 7$ be a prime. Then $7 \nmid \tau(p^2)$.
\end{lem}
\begin{proof}
Suppose $7 \mid \tau(p^2)$. Then by \eqref{mod7}
\[
p^{18}+p^9+1 \equiv 0 \mod{7}.
\]
But $p^{18}=(p^6)^3 \equiv 1$ and $p^9 \equiv p^3 \equiv \pm 1 \mod{7}$ giving
a contradiction.
\end{proof}
\begin{lem}\label{lem:mod5}
Let $p \ne 5$ be a prime. Then $5 \nmid \tau(p^2)$.
\end{lem}
\begin{proof}
Suppose $5 \mid \tau(p^2)$. Then by \eqref{mod5}
\[
(p^{71})^2+p^{71}+1 \equiv 0 \mod{5}.
\]
However, this contradicts the fact that  the congruence $x^2+x+1 \equiv 0 \mod{5}$ has no solutions.
\end{proof}
\begin{lem}\label{lem:mod9}
Let $p \ne 3$ be a prime. Then $9 \nmid \tau(p^2)$.
\end{lem}
\begin{proof}
Suppose $9 \mid \tau(p^2)$. Then by \eqref{mod3}
\[
(p^{1231})^2+p^{1231}+1 \equiv 0 \mod{9}.
\]
Since the congruence $x^2+x+1 \equiv 0 \mod{9}$ has no solutions, we obtain the desired contradiction.
\end{proof}

\section{Lucas Sequences}\label{sec:lucas}

In this section, for a fixed prime $p$ with $\tau(p) \ne 0$,
we show that the sequence $m \mapsto \tau(p^{m-1})$ can be appropriately
scaled to yield a Lucas sequence.
We begin by introducing Lucas sequences and recalling some of their
properties, 
mostly following the article of Bilu, Hanrot and Voutier \cite{BHV}.
A \textbf{Lucas pair} is a pair $(\alpha,\beta)$ of algebraic
numbers such that $\alpha+\beta$ and $\alpha\beta$ are
non-zero coprime rational integers, and $\alpha/\beta$
is not a root of unity. In particular, associated
to the Lucas pair $(\alpha,\beta)$ is a
\textbf{characteristic polynomial}
\[
X^2-(\alpha+\beta) X + \alpha \beta \; \in \; \Z[X].
\]
This polynomial has discriminant $D=(\alpha-\beta)^2 \in \Z \setminus \{0\}$.
Given a Lucas pair $(\alpha,\beta)$,
the corresponding \textbf{Lucas sequence} is given by
\[
u_n=\frac{\alpha^n-\beta^n}{\alpha-\beta}, \qquad n=0,1,2,\dotsc.
\]
Let $(\alpha,\beta)$ be a Lucas pair. A prime $\ell$ is a
\textbf{primitive divisor} of the $n$-th term of the corresponding
Lucas sequence
if $\ell$ divides $u_n$ but $\ell$ fails to divide $(\alpha-\beta)^2 \cdot u_1 u_2 \dotsc u_{n-1}$.
We shall make
essential use of the celebrated Primitive Divisor Theorem
of Bilu, Hanrot and Voutier \cite{BHV}.
\begin{thm}[Bilu, Hanrot and Voutier]\label{thm:BHV}
Let $(\alpha,\beta)$ be a Lucas pair.
If $n \ge 5$ and $n \ne 6$ then $u_n$ has a primitive divisor.
\end{thm}

Let $\ell$ be a prime. The smallest positive integer $m$
such that $\ell \mid u_m$ is called the \textbf{rank of
apparition of $\ell$}; we denote this
by $m_\ell$.
We shall also have need of  the following classical theorem of Carmichael \cite{Car}.
\begin{thm}[Carmichael]\label{thm:Carmichael}
Let $(\alpha,\beta)$ be a Lucas pair and $\ell$ be a prime.
\begin{enumerate}[label=(\roman*)]
\item If $\ell \mid \alpha \beta$ then $\ell \nmid u_m$
for all positive integers $m$.
\item Suppose $\ell \nmid \alpha \beta$.
Write $D=(\alpha-\beta)^2 \in \Z$.
\begin{enumerate}[label=(\alph*)]
\item If $\ell \ne 2$ and $\ell \mid D$, then
$m_\ell=\ell$.
\item If $\ell \ne 2$ and $\left(\frac{D}{\ell}\right)=1$, then $m_\ell \mid (\ell-1)$.
\item If $\ell \ne 2$ and $\left(\frac{D}{\ell} \right)=-1$, then $m_\ell \mid (\ell+1)$.
\item If $\ell=2$, then $m_\ell=2$ or $3$.
\end{enumerate}
\item If $\ell \nmid \alpha \beta$ then
\[
\ell \mid u_m \iff m_\ell \mid m.
\]
\end{enumerate}
\end{thm}
\begin{proof}
Note that the sequence $\{u_n\}$ satisfies the recurrence
\[
u_{n+2}-(\alpha+\beta) u_{n+1}+\alpha \beta u_n=0, \qquad u_0=0, \quad
u_1=1.
\]
	If $\ell \mid \alpha \beta$ then $u_n \equiv (\alpha+\beta)^{n-1} \pmod{\ell}$
for all $n \ge 1$. Since $\alpha+\beta$ and $\alpha \beta$
are coprime, $\ell \nmid (\alpha+\beta)$ and so $\ell \nmid u_n$
for all $n \ge 1$.

Suppose now that $\ell \nmid \alpha \beta$. Let $K=\Q(\alpha)=\Q(\beta)
=\Q(\sqrt{D})$
and $\lambda$ be a prime of $\OO_K$ above $\ell$. We first consider
(a). Here $\alpha=\beta+\gamma$ where $\lambda \mid \gamma$.
Thus
\[
u_n=\frac{\alpha^n-\beta^n}{\alpha-\beta}=
\frac{(\beta+\gamma)^n-\beta^n}{\gamma} \equiv n \beta^{n-1}  \pmod{\lambda}.
\]
Thus $\ell \mid u_n$ if and only if $\ell \mid n$.

Next we consider cases (b) and (c) together. Note that
\[
\ell \mid u_n \iff (\alpha/\beta)^n \equiv 1 \pmod{\lambda}.
\]
Thus
$m_\ell$ is equal to the multiplicative order
of the image of $\alpha/\beta$ in $(\OO_K/\lambda)^*$.
If $D$ is a quadratic residue modulo $\ell$,
then $\ell$ splits as a product of two degree $1$
primes $\lambda$, $\lambda^\prime$ of $\OO_K$.
Thus $(\OO_K/\lambda)^* \cong \F_\ell^*$ and so
$m_\ell \mid (\ell-1)$.
Finally we suppose $D$ is a quadratic non-residue modulo $\ell$.
Then $\lambda=\ell \OO_K$, and so the natural map
$\Gal(K/\Q) \rightarrow \Gal(\F_\lambda/\F_\ell)$ is an isomorphism.
Note that $\alpha$ and $\beta$ are conjugates, and so
$(\alpha/\beta)^\ell \equiv \beta/\alpha \pmod{\lambda}$.
Thus $(\alpha/\beta)^{\ell+1} \equiv 1 \pmod{\lambda}$, whence
$m_\ell \mid (\ell+1)$. The final part of the theorem
is now also clear.
\end{proof}

\bigskip

Let us fix a prime $p$ and consider
the sequence 
\begin{equation} \label{Seq}
\left\{ 1, \tau (p), \tau (p^2), \tau (p^3), \ldots \right\}.
\end{equation}
We will associate to  this a Lucas pair and a corresponding Lucas sequence.
Our starting point is the identity
\begin{equation}\label{eqn:recursion}
\tau (p^m) = \tau (p) \tau (p^{m-1}) - p^{11} \tau (p^{m-2}),
\end{equation}
valid for all integer $m \geq 2$. Once again, this was conjectured by Ramanujan
\cite{Ram} and proved by Mordell \cite{Mor}.

Let $\gamma$ and
$\delta$ be the
roots of the quadratic equation
\[
X^2- \tau (p) X + p^{11}=0,
\]
so that
\[
\gamma+\delta=\tau (p) \qquad \text{ and } \qquad
\gamma \delta = p^{11}.
\]
Then
\[
(\gamma-\delta)^2 = (\gamma+\delta)^2 - 4
\gamma \delta = \tau^2 (p) - 4 p^{11}.
\]
It follows from
Deligne's bounds (\ref{Del:bound})  that  $\gamma$
and $\delta$ are non-real Galois conjugates.
From (\ref{eqn:recursion}), we have
\begin{equation}\label{eqn:taulucas}
\tau (p^m)
=\frac{{\gamma}^{m+1}-{\delta}^{m+1}}{\gamma-\delta}.
\end{equation}
\begin{lem}\label{lem:tauval}
Suppose $\tau(p) \ne 0$. Then $\ord_p(\tau(p)) \le 5$.
\end{lem}
\begin{proof}
From (\ref{Del:bound}), 
if $p^6 \mid \tau(p)$ and $\tau(p) \ne 0$, then necessarily 
$p \le 3$. However,
\[
\tau(2)=-2^3 \times 3 \; \; \mbox{ and } \; \;  \tau(3)=2^2 \times 3^2 \times 7,
\]
providing a contradiction and completing the proof.
\end{proof}
\begin{lem}\label{lem:notrootofunity}
Suppose $\tau(p) \ne 0$. Then $\gamma/\delta$
is not a root of unity.
\end{lem}
\begin{proof}
Observe that
\[
\frac{\gamma}{\delta}+
\frac{\delta}{\gamma}+2
\; = \; \frac{\tau(p)^2}{p^{11}}.
\]
By the previous lemma, the rational number $\tau(p)^2/p^{11}$
is not an integer and therefore not an algebraic integer.
It follows that $\gamma/\delta$ is not a root
of unity.
\end{proof}
The following is an immediate
consequence of \eqref{eqn:taulucas} and Lemma~\ref{lem:notrootofunity}.
\begin{lem}\label{taup=0}
If $\tau(p)=0$ then
\[
\tau(p^{m})=\begin{cases}
0 & \text{if $m$ is odd},\\
(-p^{11})^{m/2} & \text{if $m$ is even}.
\end{cases}
\]
If $\tau(p) \ne 0$ then $\tau(p^m) \ne 0$ for all $m \ge 1$.
\end{lem}

Note that $\gcd(\gamma+ \delta,\gamma \delta)
=\gcd(\tau(p),p^{11})=1$
if and only if $p \nmid \tau(p)$. Thus
the sequence $m \mapsto \tau(p^{m-1})$ is a Lucas
sequence precisely when $p \nmid \tau(p)$. We note that $p \mid \tau(p)$
for
$$
p =  2, 3, 5, 7, 2411,  7758337633, \ldots
$$
We expect that $p \mid \tau(p)$
for infinitely many primes $p$; see Lygeros and Rozier \cite{LyRo0} for a
discussion of this problem and related computations. We will scale the pair
$(\gamma,\delta)$ to obtain a Lucas pair, 
provided $\tau(p) \ne
0$.
\begin{lem}\label{lem:associated}
Suppose $\tau(p) \ne 0$. Write $r=\ord_p(\tau(p))$ and
let
\[
\alpha=\frac{\gamma}{p^r}, \qquad \beta=\frac{\delta}{p^r}.
\]
Then $(\alpha,\beta)$ is a Lucas pair. Denoting the corresponding
Lucas sequence by $u_n$, we have
\begin{equation}\label{eqn:unform}
u_n=\frac{\tau(p^{n-1})}{p^{r(n-1)}}, \qquad n \ge 1.
\end{equation}
Moreover, $p \nmid u_n$ for all $n \ge 1$.
\end{lem}
\begin{proof}
Note that $\alpha+\beta=\tau(p)/p^{r}$ and $\alpha \beta=p^{11-2r}$
are coprime rational integers thanks to Lemma~\ref{lem:tauval}.
The identity \eqref{eqn:unform} follows immediately
from \eqref{eqn:taulucas}.
The last part is a consequence of part (i) of Theorem~\ref{thm:Carmichael}
since $p \mid \alpha \beta$.

For future reference, we note that, for $\{ u_n \}$, we have
\begin{equation} \label{disc}
D = \left( \alpha-\beta \right)^2 = p^{-2r} \left(  \tau^2 (p) - 4 p^{11} \right).
\end{equation}
\end{proof}

\section{On three sequences of polynomials}\label{sec:Fn}
We begin by defining, for $m \ge 0$, a sequence of polynomials $H_m(Z,W) \in \Z[Z,W]$
\begin{equation}\label{eqn:Hndef}
	H_m(Z,W)=\begin{cases}
		(Z^m-W^m)/(Z-W) \qquad \text{if $m$ is odd}\\
		(Z^m-W^m)/(Z^2-W^2) \qquad \text{if $m$ is even}.
		\end{cases}
\end{equation}
Let $G$ be the group generated the involutions $\kappa_1$ and $\kappa_2$ on $\Z[Z,W]$ given by
\[
\kappa_1 \; : \; 
\begin{cases}
Z \mapsto -Z, \\
W \mapsto -W,
\end{cases}
\qquad 
\kappa_2 \; : \; 
\begin{cases}
Z \mapsto W, \\ 
W \mapsto Z.
\end{cases}
\]
We compute the subring of invariants $\Z[Z,W]^G$.
\begin{lem}
$\Z[Z,W]^G=\Z[ZW,(Z+W)^2]$.
\end{lem}
\begin{proof}
It is clear that $\Z[ZW,(Z+W)^2]$ belongs to the ring of invariants.
Let $F \in \Z[Z,W]$ belong to the ring of invariants. We would
like to show that $F \in \Z[ZW,(Z+W)^2]$. Observe that $\kappa_1$ and $\kappa_2$
send monomials to monomials and preserve the degree.
Thus every homogenous component of $F$ belongs to the ring of invariants,
and we may suppose that $F$ is homogeneous. As $F$ is invariant
under $\kappa_1$ it has even degree, $2n$ say, and we may write
\[
	F=\sum_{k=0}^{2n} a_k Z^{2n} W^{2n-k}.
\]
As $F$ is invariant under $\kappa_2$ we have $a_k=a_{2n-k}$ for $k=0,\dotsc,n$.
Thus
\[
	F=a_0(Z^{2n}+W^{2n})+a_1 (ZW) (Z^{2n-2}+W^{2n-2})+a_2 (ZW)^2 (Z^{2n-4}+W^{2n-4})+\cdots+a_n (ZW)^n.
\]
	To complete the proof all we need to show is that $Z^{2n}+W^{2n} \in \Z[ZW,(Z+W)^2]$ for all $n$.
This follows from an easy induction using the identity
\[
	Z^{2n}+W^{2n} \; =\; ((Z+W)^2-2ZW) \cdot (Z^{2n-2}+W^{2n-2}) \, - \, (ZW)^2 \cdot (Z^{2n-4}+W^{2n-4}).
\]
\end{proof}
Note that the $H_m(Z,W)$ are invariant under $\kappa_1$, $\kappa_2$ and so belongs to the invariant ring $\Z[ZW,(Z+W)^2]$.
It follows that there is a sequence of polynomials $F_m(X,Y) \in \Z[X,Y]$ such that
\begin{equation}\label{eqn:FnHn}
	F_m(ZW,(Z+W)^2)=H_m(Z,W).
\end{equation}
The following lemma aids in the computation of the $F_m$.
\begin{lem}
The sequence $F_m(X,Y) \in \Z[X,Y]$ satisfies	
\[
F_0=0, \quad F_1=F_2=1, \quad F_3=-X+Y, 
\]
and the recurrence
\begin{equation} \label{eqn:Fnrecurrence}
	F_{m+2}(X,Y)=(-2X+Y)\cdot F_m(X,Y)-X^2\cdot F_{m-2}(X,Y), \qquad \text{for $m \ge 2$}. 
\end{equation}
\end{lem}
\begin{proof}
	Since $H_0=0$ and $H_1=H_2=1$ we have $F_0=0$ and $F_1=F_2=1$. Moreover, $H_3=Z^2+ZW+W^2=-ZW+(Z+W)^2$
	so $F_3=-X+Y$.
The map 
\[
	\Z[X,Y] \rightarrow \Z[ZW,(Z+W)^2], \qquad X \mapsto ZW, \quad Y \mapsto (Z+W)^2
\]
	is an isomorphism of rings that sends $F_m(X,Y)$ to $H_m(Z,W)$. Applying this isomorphism to  \eqref{eqn:Fnrecurrence}
	gives
\[
	H_{m+2}(Z,W)=(Z^2+W^2) H_m (Z,W)-(ZW)^2 H_{m-2}(Z,W)
\]
and so it is enough to prove this identity. However this identity easily follows from the definition of $H_m$
	in \eqref{eqn:Hndef}.
\end{proof}

\begin{lem}
If $m$ and $n$ are positive integers, then
\begin{enumerate}[label=(\roman*)]
	\item $F_m$ is homogeneous of degree $\lfloor (m-1)/2 \rfloor$.
	\hskip12ex \item $F_n \mid F_m$ whenever $n \mid m$.
\end{enumerate}		
\end{lem}
\begin{proof}
These follow immediately from the corresponding properties for the $H_m$.
\end{proof}

\begin{lem}\label{lem:split}
Let $m \ge 3$. Then
\begin{equation}\label{eqn:FmFactor}
	F_m(X,Y) = \prod_{j=1}^{\lfloor (m-1)/2 \rfloor} (Y-4\cos^2(\pi j/m) X)
\end{equation}
\end{lem}
\begin{proof}
Fix $m \ge 3$ and write
$\zeta=\exp(2\pi i/m)$.
Note that
\[
\begin{split}
	H_m(Z,W) &=\prod_{j=1}^{\lfloor (m-1)/2 \rfloor} (Z-\zeta^j W)(Z-\zeta^{-j} W)\\
		&=\prod_{j=1}^{\lfloor (m-1)/2 \rfloor} (Z^2+W^2-(\zeta^j+\zeta^{-j}) ZW)\\
		&=\prod_{j=1}^{\lfloor (m-1)/2 \rfloor} ((Z+W)^2-(\zeta^j+\zeta^{-j}+2) ZW)\\
		&=\prod_{j=1}^{\lfloor (m-1)/2 \rfloor} ((Z+W)^2-(2+2\cos{2\pi j/m}) ZW)\\
		&=\prod_{j=1}^{\lfloor (m-1)/2 \rfloor} ((Z+W)^2- 4 \cos^2(\pi j/m) ZW).\\
\end{split}
\]
The lemma follows.
\end{proof}

Next we define
\begin{equation}\label{eqn:PsimFactor}
	\Psi_m(X,Y)=\prod_{\stackrel{j=1}{(j,m)=1}}^{\lfloor (m-1)/2 \rfloor} (Y-4 \cos^2(\pi j/m) X).
\end{equation}
Note that $\Psi_m(X,Y) \in \Z[X,Y]$. Indeed,
\[
	\Psi_m(X,Y)=\frac{F_m(X,Y)}{\LCM\{F_n(X,Y) \: : \; n \mid m,~ n<m\}}.
\]
It follows that $\Psi_m(X,Y) \mid \Psi_n(X,Y)$ (with the divisibility being valid
in $\Z[X,Y]$) whenever $m \mid n$.
From \eqref{eqn:FmFactor} and \eqref{eqn:PsimFactor}, we see that
\[
	F_m(X,Y)=\prod_{d \mid m} \Psi_d(X,Y).
\]
We deduce that
\begin{equation}\label{eqn:inversion}
	\Psi_m(X,Y)=\prod_{d \mid m} F_d(X,Y)^{\mu(m/d)}
\end{equation}
where $\mu$ denotes the M\"{o}bius function.
\begin{lem}\label{lem:abelian}
Let $m \ge 3$ and write $\zeta_m=\exp(2\pi i/m)$.
	The polynomial $\Psi_m(1,Y)$ is monic and irreducible
	of degree $\phi(m)/2$. It is a defining polynomial
	for the abelian extension $\Q(\zeta_m)^+/\Q$.
\end{lem}
\begin{proof}
Recall that $\Q(\zeta_m)^+=\Q(\zeta_m+\zeta_m^{-1})$. 
The elements of the Galois group for $\Q(\zeta_m)/\Q$ 
are the automorphisms $\sigma_j : \zeta_m \mapsto \zeta_m^j$
with $\gcd(j,m)=1$.
Therefore Galois conjugates
of $\zeta_m+\zeta_m^{-1}+2$ are precisely $\zeta_m^j+\zeta_m^{-j}+2$ with $\gcd (j,m)=1$,
and these are the roots of $\Psi_m(1,Y)$.
The lemma follows.
\end{proof}
We shall need the following weak bound for the coefficients of $\Psi_m$.
\begin{lem}\label{lem:Hbound}
	The coefficients of $\Psi_m$ are bounded in absolute value by $5^{\phi(m)/2}$.
\end{lem}
\begin{proof}
The roots of the monic polynomial $\Psi_m(1,Y)$ are bounded in
absolute value by $4$. Writing $\Psi_m(1,Y)=\sum a_i Y^i$
	and $(4+Y)^{\phi(m)/2}=\sum b_i Y^i$
	we have $\lvert a_i \rvert \le b_i$.
Thus 
$$
\sum \lvert a_i \rvert \le \sum b_i =5^{\phi(m)/2}.
$$
\end{proof}
\begin{lem}\label{lem:hR}
Let $m \ge 3$ and write $\zeta_m=\exp(2\pi i/m)$.
Write $h_m$ and $R_m$ for the class
number and regulator of $K_m$. Then as $m \rightarrow \infty$,
\[
	\log(h_m)=O(m\log{m}), \qquad \log(h_m R_m)=O(m \log{m})
\]
where the implicit constants are absolute and effective.
\end{lem}
\begin{proof}
Write $d=\phi(m)/2$ for the degree of $K_m$.
By \cite[Propostion 2.7]{Wa} and \cite[Lemma 4.19]{Wa},
\[
	\log(\lvert \Disc(K_m)\rvert) \; \le\;
	\frac{1}{2} \log (\Disc(\Q(\zeta_m))) 
	\; \le \;
	\frac{1}{2}\phi(m) \log(m) .
\]
A theorem of Lenstra \cite[Theorem 6.5]{Lenstra} asserts that
for a number field $K$ of degree $d \ge 2$, signature $(r,s)$,
absolute discriminant $D$, class number $h$ and regulator $R$,
\[
	h \; \le \; \frac{1}{(d-1)!} \cdot \Delta \cdot 
        (d-1+\log{\Delta})^s,
        \qquad \Delta=(2/\pi)^s \cdot D^{1/2}.
\]
and
\[
	hR \; \le \; \frac{1}{(d-1)!} \cdot \Delta \cdot 
	(d-1+\log{\Delta})^s \cdot
	(\log{\Delta})^{d-1-s} . 
\]
We take $K=K_m$, so $d=\phi(m)/2$, $s=0$,
and $\Delta=\lvert \Disc(K_m) \rvert^{1/2}$.
The lemma follows.
\end{proof}
We can also deduce the bound $\log(h_m R_m)=O(m \log{m})$
from the Brauer--Siegel theorem, at the cost of introducing ineffectivity.

\begin{lem}\label{lem:tauFm}
Let $p$ be a prime. Then, for $m \ge 1$,
\[
	\tau(p^{m-1})=\tau(p)^\varepsilon \cdot F_m(p^{11},\tau(p)^2), \qquad
	\varepsilon=\begin{cases}
		0 & \text{if $m$ is odd}\\
		1 & \text{if $m$ is even}.
	\end{cases}
\]
In particular, $\Psi_m(p^{11},\tau(p)^2) \mid \tau(p^{m-1})$.
\end{lem}
\begin{proof}
From \eqref{eqn:taulucas}
\[
\tau(p^{m-1})=
\frac{{\gamma}^{m}-{\delta}^m}{\gamma-\delta}=\begin{cases}
	H_m(\gamma,\delta) & \text{if $m$ is odd}\\
	(\gamma+\delta) H_m(\gamma,\delta) & \text{if $m$ is even}, 
\end{cases}
\]
where $\gamma+\delta=\tau(p)$ and $\gamma\delta= p^{11}$.
The lemma follows from \eqref{eqn:FnHn}.
\end{proof}

\medskip

\begin{lem}\label{lem:ramification}
Let $m=5$ or $m \ge 7$. Then precisely the same primes
ramify in $L_m=\Q(\zeta_m)$
as in $K_m=\Q(\zeta_m)^+$.
\end{lem}
\begin{proof}
Note that for $m \in \{2, 3, 4, 6\}$, we have $K_m=\Q$
so the conclusion of the lemma is false in those cases.

By the proof of Proposition 2.15 of \cite{Wa}, we know that
$L_m/K_m$ is unramified if $m$ is divisible by at least two distinct
odd primes, or divisible by $4$ and an odd prime. 
We may therefore suppose that $m \in \{ 2^a, p^a, 2 p^a\}$ where
$p$ is an odd prime and $a$ is a positive integer. If $m=2^a$ with $a \ge 3$, then
$K_m$ has degree $2^{a-2}>1$, and the set of ramified
primes for both $L_m$ and $K_m$ is $\{2\}$. 
Let $p$ be an odd prime and $a \ge 1$.
Then $L_{2p^{a}}=L_{p^{a}}$ and $K_{2 p^a}=K_{p^a}$.
Now the set of ramified primes for $L_{p^{a}}$
and $K_{p^a}$ is just $\{p\}$ as long
as the degree $\phi(p^a)/2$ of $K_{p^{a}}$
is $>1$. The lemma follows.
\end{proof}

\begin{lem}\label{lem:preThue}
Let $m =5$ or $m \ge 7$. Let $x$ and $y$ be coprime integers, and $q$
be a prime. Suppose $q^a \mid\mid \Psi_m(x,y)$ with $a \ge 1$ an integer. Then
either $q \equiv \pm 1 \pmod{m}$ or $q^a \mid m$. 
\end{lem}
\begin{proof}
Write $L_m=\Q(\zeta_m)$. Recall the isomorphism
\[
	(\Z/m\Z)^* \cong \Gal(L_m/\Q), \qquad j \; \mapsto\;  
	(\sigma_j \; : \; \zeta_m \mapsto \zeta_m^{j}).
\]
The subfield $K_m=\Q(\zeta_m+\zeta_m^{-1})$ is the fixed
	field for $\langle \sigma_{-1} \rangle=\{\sigma_1,\sigma_{-1}\}$.
Let $q$ be a rational prime that does not ramify in $K_m$
(and hence in $L_m$ by Lemma~\ref{lem:ramification}). 
The Frobenius automorphism
for $q$ is simply $\sigma_q$. The prime $q$ splits completely
in $K_m$ if and only if the restriction of $\sigma_q$ to $K_m$
is trivial. This is equivalent to 
	$\sigma_q \in \{\sigma_1,\sigma_{-1}\}$
and therefore equivalent to
$q \equiv \pm 1 \pmod{m}$.

Let $\lambda=\zeta_m+\zeta_m^{-1}+2$. This is a root
for $\Psi_m(1,Y)$ and also a generator for $K_m$. 
Note \cite[Proposition 2.16]{Wa} that $\OO_{K_m}=\Z[\lambda]$.
We are given that $q \mid \Psi_m(x,y)$.
Since $\Psi_m(1,Y)$ is monic, if $q \mid x$ then $q \mid y$
giving a contradiction. Hence $q \nmid x$ and so
$\Psi_m(1,y/x) \equiv 0 \pmod{q}$.
By the Dedekind-Kummer Theorem, there is a degree $1$ prime
ideal $\fq$ above $q$. As $K_m/\Q$ is Galois, all primes
above $q$ must therefore have degree $1$. Thus $q$ is either totally
split or ramified in $K_m$. 
If $q$ is totally split, then $q \equiv \pm 1 \pmod{m}$
and we are finished.

We shall therefore suppose that $q$ is ramified in $K_m$.
Let $\fq_1,\dotsc,\fq_r$ be the prime ideals of $\OO_{K_m}$
above $q$. Write $G$ for the Galois group of $\Gal(K_m/Q)$,
and let $I$ be the inertia subgroup for $\fq_1$. 
As $G$ is abelian, $I$ is also the inertia subgroup
for all $\fq_i$. Thus $\fq_i^\sigma=\fq_i$
for all $\sigma \in I$ and for $i=1,\dots,r$.
	Since $q$ is ramified, $I \ne 1$. Fix $\sigma_j \in \Gal(L_m/\Q)$
whose restriction to $K_m$ is a non-trivial element of $I$.
Thus $\gcd (j,m)=1$ and $j \not \equiv \pm 1 \pmod{m}$.
Write $\lambda_j=\zeta_m^{j}+\zeta_m^{-j}+2=\sigma_j(\lambda)$.

We factor the ideal $(y-\lambda x)\OO_{K_m}$
as
\begin{equation}\label{eqn:fafb}
	(y-\lambda x)\OO_{K_m}=\fa \fb
\end{equation}
where $\fa$, $\fb$ are ideals with $\fa$ supported
on $\fq_1,\dotsc,\fq_r$, and $\fb$ not divisible
by $\fq_1,\dotsc,\fq_r$. 
By assumption $q^a \mid\mid \Psi_m(x,y)$.
However $\Psi_m(x,y)=\Norm_{K_m/\Q}(y-\lambda x)$
and thus $\Norm_{K_m/\Q}(\fa)=q^a$.
Note that any ideal dividing both $x$ and $\fa$
must also divide $y$ by \eqref{eqn:fafb}.
As $x$, $y$ are coprime, we deduce that $x$ and $\fa$
are coprime.

Since $\fq^{\sigma_j}=\fq$ for all $\fq \mid \fa$, 
we have $\fa^{\sigma_j}=\fa$.  
Hence $\fa$ divides
\[
	(y-\lambda^{\sigma_j} x) - (y - \lambda x)= (\lambda-\lambda^{\sigma_j})x. 
\]
Thus $\fa$ divides
\[
	\lambda-\lambda^{\sigma_j}=
	(\zeta_m+\zeta_m^{-1})-(\zeta_m^{j}+\zeta_m^{-j})
	=\zeta_m^{-1} (\zeta_m^{j+1}-1) (\zeta_m^{-j+1}-1)
\]
and it follows that $q^{2a}=\Norm_{K_m/K}(\fa)^2$ divides
\[
	\begin{split}
		(\Norm_{K_m/\Q}(\lambda-\lambda^{\sigma_j}))^2 &=
	\Norm_{L_m/\Q}(\lambda-\lambda^{\sigma_j})\\
		&=
		\Norm_{L_m/\Q}(\zeta_m^{j+1}-1) \cdot
		\Norm_{L_m/\Q}(\zeta_m^{-j+1}-1). \\
	\end{split}
\]
This divides $\prod_{i=1}^{m-1} (\zeta_m^i-1)^2=m^2$.
Hence $q^a \mid m$ as required.
\end{proof}

\begin{lem}\label{lem:psiDiv}
Let $p$ be a prime and suppose $\tau(p)=0$.
Let $r=\ord_p(\tau(p))$ and write
\[
	x=p^{11-2r} \; \mbox{ and } \;  y=\frac{\tau(p)^2}{p^{2r}}.
\]
Let $\{u_m\}$ be the Lucas sequence defined in Lemma~\ref{lem:associated}. Then,
for $m \ge 3$,
\begin{equation}\label{eqn:psiu}
	\Psi_m(x,y)= \prod_{d \mid m} u_d^{\mu(m/d)}.
\end{equation}
Moreover, if $m =5$ or $m \ge 7$, then $\Psi_m(x,y)$
is divisible by some prime $\ell \nmid m$.
\end{lem}
\begin{proof}
Note that $x$ and $y$ are in fact integers by Lemma~\ref{lem:tauval},
and are coprime by the definition of $r$.
Let $\{u_m\}$ be the Lucas sequence defined in Lemma~\ref{lem:associated}.
Thus
\[
u_m=\frac{\tau(p^{m-1})}{p^{r(m-1)}}=\frac{\alpha^m-\beta^m}{\alpha-\beta}, 
\qquad \alpha \beta=p^{11-2r}, \qquad \alpha+\beta=\tau(p)/p^r.
\]
Write
\[
\varepsilon(m)=\begin{cases}
	0 & \text{if $m$ is odd}\\
	1 & \text{if $m$ is even}.
\end{cases}
\]
Then
\[
\begin{split}
	F_m(x,y) &=\frac{1}{p^{2 r \deg(F_m)}} \cdot F_m(p^{11}, \tau(p)^2)\\
		& =\frac{\tau(p^{m-1})}{p^{2 r \deg(F_m) } \cdot \tau(p)^{\varepsilon(m)}}
	\qquad \text{(by Lemma~\ref{lem:tauFm})}
		\\ 
		&=\frac{p^{r(m-1)} \cdot u_m}{p^{2 r \deg(F_m) } \cdot \tau(p)^{\varepsilon(m)}} .
	\end{split}
\]
However, since $\deg(F_m)=\lfloor (m-1)/2 \rfloor$, it follows that
\[
	F_m(x,y)=\left( \frac{p^r}{\tau(p)} \right)^{\varepsilon(m)} \cdot u_m.
\]
By \eqref{eqn:inversion},
\[
	\Psi_m(x,y)=\left( \frac{p^r}{\tau(p)} \right)^{\sum_{d \mid m} \varepsilon(d) \mu(m/d)} \cdot \prod_{d \mid m} u_d^{\mu(m/d)}.
\]
It is easy to see that 
\[
	\sum_{d \mid m} \varepsilon(d) \mu(m/d)=\begin{cases}
		0 & \text{if $m \ne 2$}\\
		1 & \text{if $m=2$}.
	\end{cases}
\]
This completes the proof of \eqref{eqn:psiu}.

Now let $m=5$ or $m \ge 7$. By Theorem~\ref{thm:BHV}, the term $u_m$
has a prime divisor $\ell$ that does not divide $(\alpha-\beta)^2$
nor $u_1 u_2\cdots u_{m-1}$. 
By Theorem~\ref{thm:Carmichael},
we know that $\ell \ne p$, that $m=m_\ell$,
and that $m\mid (\ell-1)$
or $m \mid (\ell+1)$. In particular, $\ell \nmid m$.
From \eqref{eqn:psiu}, we have $\ell \mid \Psi_m(x,y)$
as required.
\end{proof}

\section{Proof of Theorem~\ref{thm:main2}}\label{sec:Bugeaud}
We shall need the following theorem \cite[Theorem 1]{Bug-BLMS}.
\begin{thm}[Bugeaud]
Let $K$ be a number field. Let $u \ge 2$ and $v \ge 3$ be integers, and 
let $a$, $b\in \OO_K\setminus \{0\}$.
There exist effectively computable positive constants $\varepsilon_1$, $\varepsilon_2$ depending only on $a$, $b$, $u$, $v$ and $K$
such that every pair of coprime $x$, $y \in \OO_K$ with 
\[
	\max\{\lvert \Norm_{K/\Q}(x) \rvert,~\lvert \Norm_{K/\Q}(y) \rvert\} \; > \; \varepsilon_1
\]
satisfy
\[
	P(a x^u+by^v) \; > \; \varepsilon_2 \cdot \log\log \max\{ \lvert \Norm_{K/\Q}(x)\rvert,~\lvert \Norm_{K/\Q}(y) \rvert\}.
\]
\end{thm}
In the above theorem, $P(\delta)$ for $\delta \in \OO_K$ denotes the largest rational prime
that is below a prime ideal dividing $\delta$.

\bigskip

We now prove Theorem~\ref{thm:main2}.
Let $p$ be a prime and suppose $\tau(p) \ne 0$.
Let $m \ge 3$. We want to show that
\[
	P(\tau(p^{m-1})) \; \gg_{m} \; \log\log{p}.
\]
Let $r=\ord_p(\tau(p))$ and recall that $\ord_p(\tau(p^{m-1}))=r(m-1)$ by Lemma~\ref{lem:tauval}.
If $r \ge 1$ then
\[
	P(\tau(p^{m-1})) \; \ge \; p,
\]
whereby we may suppose that $r=0$. 
Recall that $\Psi_m(p^{11},\tau(p)^2) \mid \tau(p^{m-1})$ by Lemma~\ref{lem:tauFm}.
Let $K=K_m=\Q(\zeta_m)^+$ and let $\lambda=\zeta_m+\zeta_m^{-1}+2$
which is a root of the monic polynomial $\Psi_m(1,Y)$. 
Then 
\[
	\Psi_m(p^{11},\tau(p)^2) \; =\; \Norm_{K/\Q}(\tau(p)^2- \lambda \cdot p^{11})
\]
and therefore 
\[
	P(\tau(p^{m-1})) \; \ge \; P(\Psi_m(p^{11},\tau(p)^2) ) \; = \; 
	P(\tau(p)^2-\lambda \cdot p^{11}).
\]
We now apply Bugeaud's theorem with $u=2$, $v=11$, $a=1$, $b=-\lambda$,
$x=\tau(p)$, $y=p$ to deduce that $P(\tau(p)^2- \lambda \cdot p^{11}) \; \gg_m \; \log\log{p}$. 
This completes the proof.
\section{Proof of Theorem \ref{thm:main}}\label{sec:thm1a}
In this section, we prove Theorem~\ref{thm:main}.
For this we appeal to 
a result of Bugeaud and Gy\H{o}ry \cite{BuGy}
which provides bounds for solutions to Thue--Mahler
equations.
Let $F(X,Y) \in \Z[X,Y]$ be an irreducible binary form of degree $n \ge 3$,
and let $b$ a non-zero rational integer with absolute value
at most $B\ge e$. 
Let $H \ge 3$ be an upper bound for the absolute values
of the coefficients of $F$. 

Let $\alpha_1$, $\alpha_2$ and $\alpha_3$ be three distinct roots of $F(1,Y)$.
Define
\[
	\mathbb{M}=\Q(\alpha_1), \; \;  \mathbb{M}_{123}=\Q(\alpha_1,\alpha_2,\alpha_3) \; \;
	\mbox{ and } \; \; N=[\mathbb{M}_{123} : \Q].
\]
Write $h_{\mathbb{M}}$ for the class number of $\mathbb{M}$ and $R_{\mathbb{M}}$ for its regulator.
Let $p_1,p_2,\dotsc,p_s$ ($s>0$) be distinct primes not exceeding $P$. 
Consider the \textbf{Thue--Mahler equation}
\begin{equation}\label{eqn:TM}
	F(x,y) \; =\; b \cdot p_1^{z_1} p_2^{z_2} \cdots p_s^{z_s}, \qquad x,~y,~z_i \in \Z, \qquad \gcd(x,y,p_1p_2\cdots p_s)=1.
\end{equation}
For a positive real number $a$, we write $\log^*{a}=\max\{1,\log{a}\}$.
\begin{thm}[Bugeaud and Gy\H{o}ry]\label{thm:BG}
All solutions to \eqref{eqn:TM} satisfy
\begin{multline*}
	\log \max\{ \lvert x \rvert,~\lvert y \rvert,~ p_1^{z_1} \cdots p_k^{z_k}\} \; \le  
	\\
	c(n,s) \cdot P^N \cdot (\log{P})^{ns+2} \cdot R_{\mathbb{M}} h_{\mathbb{M}} \cdot (\log^*(R_{\mathbb{M}} h_{\mathbb{M}}))^2
	\cdot (R_\mathbb{M}+s h_\mathbb{M}+\log(HB)),
\end{multline*}
where
\[
	c(n,s)=3^{n(2s+1)+27} \cdot n^{2n (7s+13)+13} \cdot (s+1)^{5n(s+1)+15}.
\]
\end{thm}
The theorem as stated is the first part of Theorem 4 in \cite{BuGy}, with only one
minor difference. In \cite{BuGy} the authors take $N=n(n-1)(n-2)$.
However, in their proof $N$ is simply taken as an upper bound for the 
degree $[\mathbb{M}_{123} : \Q]$, and so we can take $N=[\mathbb{M}_{123}:\Q]$.

\bigskip

We now embark on the proof of Theorem~\ref{thm:main}.
In what follows $\eta_2,\eta_3,\dots$ will denote absolute
effectively computable positive constants.
Let us fix a prime $p$ and suppose $\tau(p) \neq 0$.
We will in fact show that
\begin{equation}\label{eqn:show}
	P(\tau(p^{m-1})) \; \ge \; \eta_2 \cdot \frac{\log\log(p^m)}{\log\log\log(p^m)},
\end{equation}
for $m \ge 3$ which implies \eqref{eqn:required}.
In view of Theorem~\ref{thm:main2} (which was proved in Section~\ref{sec:Bugeaud}),
we shall suppose that $m=7$, $9$, $11$ or $m \ge 13$.
In particular, $\Psi_m(X,Y)$ is irreducible
of degree $\phi(m)/2 \ge 3$. 
Write $r=\ord_p(\tau(p))$. 
By Lemma~\ref{lem:associated}, we have $\ord_p(\tau(p^{m-1}))=r(m-1)$.
Recall that $r=\ord_p(\tau(p)) \le 5$ by Lemma~\ref{lem:tauval}.
Let 
\[
	x=p^{11-2r}, \qquad y=\tau(p)^2/p^{2r},
\]
and observe that $\gcd(x,y)=1$.
By Lemma~\ref{lem:tauFm}, we know that $\Psi_m(x,y)$
is a divisor of $\tau(p^{m-1})$ and therefore
\[
	P(\tau(p^{m-1})) \; \ge \; P(\Psi_m(x,y)) 
\]
To prove \eqref{eqn:show}, we shall show that
\begin{equation}\label{eqn:desired}
	P(\Psi_m(x,y)) \; > \; \eta_3 \cdot \frac{\log\log(p^m)}{\log\log\log(p^m)}.
\end{equation}
By Lemma~\ref{lem:preThue}, we can write
\begin{equation}\label{eqn:PsiTM}
	\Psi_m(x,y)\; =\; b \cdot p_1^{z_1} p_2^{z_2} \cdots p_s^{z_s}, 
\end{equation}
where the $p_i$ are primes, and
\[
	b \mid m, \; \;  p_i \equiv \pm 1 \pmod{m} \; \; \mbox{ and } \; \; p_1<p_2<\cdots<p_s.
\]
From Lemma~\ref{lem:psiDiv}, we have $s \ge 1$.
It is clear that 
\[
	s \; < \; \eta_4 \cdot \frac{p_s}{m}.
\]
In what follows we make use of the following inequalities
\[
n \; < \; m \; \; \mbox{ and } \; \; 	ns \; < \; ms \; < \; \eta_4 \cdot p_s.
\]
Moreover, since $p_s \equiv \pm 1 \pmod{m}$, we have
\[
	p_s \; \ge \;  m-1.
\]
We will apply Theorem~\ref{thm:BG} to \eqref{eqn:PsiTM}. We take
\[
	F=\Psi_m, \; \;  B=m, \; \;  P=p_s, \; \;  n=N=\phi(m)/2 \; \; \mbox{ and } \; \; H=5^{n/2}.
\]
where the choice of $H$ is justified by Lemma~\ref{lem:Hbound}.
By Lemma~\ref{lem:hR},
we have 
\[
	\log(h_{\mathbb{M}}) \; < \; \eta_5 \cdot m \log{m} \; \; \mbox{ and } \; \; 	\log(h_{\mathbb{M}} R_{\mathbb{M}}) \; < \; \eta_6 \cdot m \log{m}.
\]
Since $x=p^{11-2r}$ with $r \le 5$, 
we have $\log\log{p} \le \log\log{x}$.
Taking logarithms in Theorem~\ref{thm:BG}, and making repeated use of the above
inequalities and bounds, yields
\[
	\log\log{p}  \; < \; \eta_7 \cdot p_s \cdot \log{p_s}. 
\]
But
\[
	\log\log(p^m) \; =
	\;
	\log\log{p}+\log{m}
	\; <\; \eta_7 \cdot p_s \cdot \log{p_s}+\log(p_s+1) \; < \; \eta_8 \cdot p_s \cdot \log{p_s}.
\]
The desired inequality \eqref{eqn:desired} follows, completing the proof 
of Theorem~\ref{thm:main}.
\section{The equation $\tau(p^2)=\kappa \cdot q^b$}\label{sec:taup2}

In this section we establish the following two propositions.
\begin{prop}\label{prop:pmqb}
Let $3 \le q <100$ be a prime.  The equation
\begin{equation}\label{eqn:pmqb}
\tau(p^2) \; = \; \pm q^b,
\qquad
\text{$p$ prime,}
\qquad \text{$b \ge 0$}
\end{equation}
has no solutions.
\end{prop}
\begin{prop}\label{prop:tpsmooth}
The equation
\begin{equation}\label{eqn:tpsmooth}
\tau(p^2) \; = \; \pm 3^{b_1} 5^{b_2} 7^{b_3} 11^{b_4},
\qquad
\text{$p$ prime,}
\qquad \text{$b_1,~b_2,~b_3,~b_4 \ge 0$}
\end{equation}
has no solutions.
\end{prop}

We consider first the following general equation.
\begin{equation}\label{eqn:quad}
\tau(p^2) \; = \; \kappa \cdot q^b,
\qquad
\text{$p \nmid 2\kappa q$ prime,}
\qquad \text{$b \ge 0$}.
\end{equation}
Here $\kappa$ is an odd integer, $q$ is an odd prime, and we assume for convenience
that $q \nmid \kappa$.
Recall that $\tau(p^2)=\tau(p)^2-p^{11}$.
Equation \eqref{eqn:quad} can be written as
\[
p^{11}+ (\kappa \cdot q^b) \cdot 1^{11}\; =\; \tau(p)^2
\]
and so is an equation of signature $(11,11,2)$.
Following the first author and Skinner \cite{BenS}, we associate to a solution of
\eqref{eqn:quad} the Frey-Hellegouarch curve
\[
\begin{cases}
E_{p}\; : \;
Y^2=X(X^2+2\tau(p) X + \tau(p)^2-p^{11}) & \text{if $p \equiv 1 \mod{4}$},\\
E_{p} \; : \;
Y^2=X(X^2+2 \tau(p) X +p^{11}) & \text{if $p \equiv 3 \mod{4}$}.\\
\end{cases}
\]
Let
\begin{equation}\label{eqn:cond}
N=\begin{cases}
2^5 \cdot \Rad(\kappa) \cdot q \cdot p & \text{if $b >0$}\\
2^5 \cdot \Rad(\kappa) \cdot p & \text{if $b=0$},
\end{cases}
\qquad
N^\prime=\begin{cases}
2^5 \cdot \Rad(\kappa) \cdot q & \text{if $11 \nmid b$}\\
2^5 \cdot \Rad(\kappa) & \text{if $11 \mid b$}.
\end{cases}
\end{equation}
Here $\Rad(\kappa)$ denotes the product of the prime divisors of $\kappa$.
The Frey-Hellegouarch curve $E_p$ has conductor $N$.
Moreover, it follows from the recipes of the first author and Skinner \cite{BenS}
(based on the modularity theorem and Ribet's level lowering theorem)
that there is a normalized newform
\begin{equation}\label{eqn:newform}
f=q+\sum_{n=1}^\infty c_n q^n
\end{equation}
of weight $2$ and level $N^\prime$
and a prime $\varpi \mid 11$ in the integers of $K=\Q(c_1,c_2,\dotsc)$
so that
\begin{equation}\label{eqn:mod11}
\overline{\rho}_{E_p,11} \sim \overline{\rho}_{f,\varpi}.
\end{equation}
The restrictions on $\kappa$ and $q$ being coprime
odd integers merely reduce the number of possibilities for $N$, $N^\prime$,
yet cover all the cases we are interested in solving.
The restriction $p \nmid 2 \kappa q$ is needed
so that the minimal discriminant $\Delta$ of the Frey-Hellegouarch curve $E_p$
satisfies $\ord_p(\Delta) \equiv 0 \mod{11}$ which is necessary for
 application of  Ribet's level lowering theorem
in order to obtain a weight $2$ newform $f$ of level $N^\prime$ not divisible
by $p$.

Throughout what follows, $\ell$ will be a prime satisfying
\begin{equation}\label{eqn:lassump}
\ell \nmid 2 \cdot 11 \cdot \kappa q  p.
\end{equation}
Then, taking traces of the images of the Frobenius element at $\ell$
in \eqref{eqn:mod11} we obtain $a_\ell(E_p) \equiv c_\ell \mod{\varpi}$ and
so
\begin{equation}\label{eqn:mod11norm}
\Norm_{K/\Q}(a_\ell(E_p)-c_\ell) \equiv 0 \mod{11}.
\end{equation}

We shall use both the congruences for the $\tau$-function
\eqref{mod2}--\eqref{mod23} and also \eqref{eqn:mod11norm}
to derive congruences for $b$.

\begin{lem}\label{lem:presieve}
Let $(p,b)$ be a solution to \eqref{eqn:quad} and suppose $p \ne 3$, $23$.
Let $\ell$ be a prime satisfying \eqref{eqn:lassump}. Let
\[
\cA_{\ell}=\{(s,t) \; : \; s,~t \in \F_\ell, \quad s \not \equiv 0 \mod{\ell},
\qquad t^2-s^{11} \not \equiv 0 \mod{\ell}
\},
\]
and
\begin{equation}\label{eqn:cB}
\cB_{\ell}=
\begin{cases}
\cA_{\ell} & \ell \ne 3,~5,~7,~23;\\
\{(s,t) \in \cA_3 \; : \; t \equiv s+1 \mod{3}\} & \ell=3;\\
\{(s,t) \in \cA_{5} \; : \; t \equiv s^2(s^{3}+1) \mod{5} \} & \ell=5;\\
\{(s,t) \in \cA_{7} \; : \; t \equiv s(s^{3}+1) \mod{7} \} & \ell=7;\\
\{(s,0) \; : \; s \in \F_{23}^*\setminus (\F_{23}^*)^2\} \cup
\{(s,t) \; : \; s \in (\F_{23}^*)^2,~t=2,~-1\}  & \ell=23.
\end{cases}
\end{equation}
For $(s,t) \in \cB_\ell$ let
\[
E_{s,t,1}/\F_\ell \; : \; Y^2=X(X^2+2 t X +t^2-s^{11}), 
\qquad
E_{s,t,3}/\F_\ell \; : \; Y^2=X(X^2+2 t X +s^{11}).
\]
Let $f$ be a newform of weight $2$ and level $N^\prime$
so that \eqref{eqn:mod11} is satisfied, and $c_\ell$
be its $\ell$-th coefficient.
For $j=1$, $3$, let
\begin{equation}\label{eqn:cC}
\cC_{\ell,j}(f)=
\{(s,t) \in \cB_\ell \; : \; \Norm(a_\ell(E_{s,t,j})-c_\ell) \equiv 0
\mod{11}\},
\end{equation}
and
\[
\cD_{\ell,j}(f)=\{ t^2-s^{11} \; : \; (s,t) \in \cC_{\ell,j}(f)\} \subseteq \F_\ell.
\]
If $p \equiv 1 \mod{4}$ then
$(\kappa \cdot q^b \mod \ell) \in \cD_{\ell,1}(f)$.
If $p \equiv 3 \mod{4}$ then
$(\kappa \cdot q^b \mod \ell) \in \cD_{\ell,3}(f)$.
\end{lem}
\begin{proof}
Since $\ell \nmid 2 \kappa q p$, and
$\tau(p)^2-p^{11}=\tau(p^2)=\kappa \cdot q^b$
we see that there is some $(s,t) \in \cA_\ell$
so that $(p,\tau(p)) \equiv (s,t) \mod{\ell}$.
Moreover, from the congruences
for $\tau$ in \eqref{mod3}--\eqref{mod23}
there is some $(s,t) \in \cB_\ell$
so that $(p,\tau(p)) \equiv (s,t) \mod{\ell}$;
it is here that we make use of the assumption $p \ne 3$, $23$.
For such a pair $(s,t)$, the
reduction modulo $\ell$ of the Frey-Hellegouarch curve $E_p$
is $E_{s,t,j}/\F_\ell$,
where $j=1$ or $3$ according to whether $p \equiv 1$ or $3 \mod{4}$.
Thus $a_\ell(E_{s,t,j}) =a_\ell(E_p)$. Hence,
$\Norm(a_{\ell}(E_{s,t,j})-c_\ell) \equiv 0 \mod{11}$
by \eqref{eqn:mod11norm}, and so $(s,t) \in \cC_{\ell,j}(f)$.
Since $t^2-s^{11} \equiv \tau(p)^2-p^{11} \equiv \kappa \cdot q^{b} \mod{\ell}$
we see that
$(\kappa \cdot q^b \bmod \ell) \in \cD_{\ell,j}(f)$.
\end{proof}

For any prime $\ell$ satisfying \eqref{eqn:lassump},
the lemma gives congruences
for $q^b$ modulo $\ell$,
and hence leads to congruences for $b$ modulo $O_{\ell}(q)$,
where $O_{\ell}(q)$ will be our notation for the multiplicative order
of $q$ modulo $\ell$. This idea is formalized
in the following lemma.
\begin{lem}\label{lem:sieve}
Let $(p,b)$ be a solution to \eqref{eqn:quad} with $p \ne 3$, $23$,
and let $M$ be a positive integer satisfying $22 \mid M$.
Define $\cE_1$ and $\cE_3$ via
\[
\cE_1=
\{ 0 \le \beta \le M-1 \; : \; \kappa \cdot q^\beta \equiv 3 \mod{4}\} 
\; \; \mbox{ and } \; \; 
\cE_3=\{ 0 \le \beta \le M-1 \; : \; \kappa \cdot q^\beta \equiv 1 \mod{4}\}.
\]
Let $f$ be a newform of weight $2$ and level $N^\prime$
so that \eqref{eqn:mod11} is satisfied.
For $j=1$, $3$, define
\[
\cF_{j}(f)=
\begin{cases}
\{ \beta \in \cE_j \; : \; 11 \nmid \beta\}
& \text{if $N^\prime=2^5 \cdot
\Rad(\kappa) \cdot q$}\\
\{ \beta \in \cE_j \; : \; 11 \mid \beta\}
& \text{if $N^\prime=2^5 \cdot
\Rad(\kappa)$}.\\
\end{cases}
\]
Suppose now that $\cL$ is a set of
primes satisfying
\begin{equation}\label{eqn:cL}
\ell \nmid 2 \cdot 11 \cdot \kappa q p, \qquad
O_{\ell}(q) \mid M.
\end{equation}
For $\ell \in \cL$ and $j=1$, $3$, let
\[
\cG_{\ell,j}(f) \; = \; \{ \beta \in \cF_{j}(f) \; : \;
(\kappa \cdot q^\beta \bmod \ell) \in \cD_{\ell,j}(f)
\}.
\]
Let
\[
\cH_j(f) \; = \; \bigcap_{\ell \in \cL} \cG_{\ell,j}(f).
\]
If $p \equiv 1 \mod{4}$ then
there is some $\beta \in \cH_1(f)$ such that
$b \equiv \beta \mod{M}$.
If $p \equiv 3 \mod{4}$ then
there is some $\beta \in \cH_3(f)$ such that
$b \equiv \beta \mod{M}$.
\end{lem}
\begin{proof}
Let $0 \le \beta \le M-1$ be the unique integer
such that $\beta \equiv b \mod{M}$.
Let $j=1$, $3$ according to whether $p \equiv 1$ or $3 \mod{4}$
respectively.
As $2 \mid M$ and $q$ is odd we have $\kappa \cdot q^\beta \equiv \kappa
\cdot q^b \mod{4}$.
Note from \eqref{mod2} that
\[
\kappa \cdot q^\beta \equiv
\kappa \cdot q^b=\tau(p^2) = \tau(p)^2-p^{11} \equiv  p^2+p+1 \equiv
\begin{cases}
3 \mod{4} & \text{if $p \equiv 1 \mod{4}$}\\
1 \mod{4} & \text{if $p \equiv 3 \mod{4}$}.\\
\end{cases}
\]
Thus $\beta \in \cE_j$.

Also $11 \mid M$. Hence $11 \mid b$
if and only if $11 \mid \beta$.
From the definition of $N^\prime$
in \eqref{eqn:cond} we see that $\beta \in \cF_j(f)$.

Now let $\ell \in \cL$. By Lemma~\ref{lem:presieve},
we know that $(\kappa \cdot q^b \bmod{\ell})
\in \cD_{\ell,j}(f)$. However $O_\ell(q) \mid M$ and $M \mid (\beta-b)$.
Thus $\kappa \cdot q^{\beta} \equiv \kappa \cdot q^{b} \mod{\ell}$,
and so $(\kappa \cdot q^\beta \bmod{\ell}) \in \cD_{\ell,j}(f)$.
We deduce that $\beta \in \cG_{\ell,j}(f)$ for all $\ell \in \cL$.
Therefore $\beta \in \cH_j(f)$ completing the proof.
\end{proof}

\begin{proof}[Proof of Proposition~\ref{prop:pmqb}]
We checked that \eqref{eqn:pmqb} has no solutions
with $p<200$ for primes $3 \le q <100$. We shall henceforth suppose
that $p>200$. In particular, $p \ne q$.
Moreover, any solution to \eqref{eqn:pmqb}
is a solution to \eqref{eqn:quad} with $\kappa=1$ or $-1$.
For a given $3 \le q < 100$ we shall let
\[
M=396=2^2 \cdot 3^2 \cdot 11, \qquad
\cL=\{\text{$3 \le \ell <200$ prime, $\ell \ne 11$, $q$} \; : \;  O_\ell(q) \mid M\}.
\]
Observe that since $p>200$ that every $\ell \in \cL$ satisfies
 \eqref{eqn:cL}.

Suppose first that $11 \mid b$ 
and write $b=11c$. Then $(x,y,z)=(p,\pm q^{c},\tau(p))$
is a solution to the equation $x^{11}+y^{11}=z^2$ satisfying $\gcd(x,y,z)=1$.
Darmon and Merel \cite{DM} showed that the equation $x^n+y^n=z^2$
has no solutions $(x,y,z) \in \Z^3$ with $n \ge 4$, $\gcd(x,y,z)=1$.
This contradiction completes the proof for $q \ne 3$. 

Thus $11 \nmid b$, and so
 in \eqref{eqn:cond} the level is $N^\prime=2^5 q$.
We will consider the case $q=3$ a little later. Suppose $5 \le q < 100$.
We wrote a \texttt{Magma} script which for each prime
$5 \le q <100$, computes the weight $2$ newforms $f$ of level
$N^\prime=2^5 q$,
and the sets $\cH_1(f)$ and $\cH_3(f)$ both for $\kappa=1$, $\kappa=-1$.
We found all of these
to be empty. By Lemma~\ref{lem:sieve}, we conclude that
\eqref{eqn:pmqb} has no solutions with $5 \le q< 100$.

\bigskip

It remains to consider the case $q=3$. By Lemma~\ref{lem:mod9}, we see that
	$b=0$ or $1$. But $11 \nmid b$, therefore $b=1$.
Thus
\begin{equation}\label{eqn:pm3}
\tau(p)^2-p^{11}=\pm 3.
\end{equation}
We consider this modulo $23$ using \eqref{mod23}.
If $p$ is a quadratic non-residue
modulo $23$, then $p^{11} \equiv -1 \mod{23}$ and $\tau(p) \equiv 0
\mod{23}$ giving a contradiction.
If $p$ is a quadratic residue
modulo $23$, then $p^{11} \equiv 1 \mod{23}$ and $\tau(p) \equiv 2$,
$-1
\mod{23}$. We conclude that $\tau(p) \equiv 2 \mod{23}$ and
$\tau(p)^2-p^{11}=3$.
Thus
\[
(\tau(p)+\sqrt{3})(\tau(p)-\sqrt{3})=p^{11}.
\]
The two factors on the left-hand side are coprime integers in
$\Z[\sqrt{3}]$. We see that
\[
\tau(p)+\sqrt{3}=(2+\sqrt{3})^a \gamma^{11}, \qquad \gamma \in \Z[\sqrt{3}],
\qquad \Norm(\gamma)=p, \qquad 0 \le a \le 10.
\]
Let $\fq=(2+3\sqrt{3})\Z[\sqrt{3}]$. Then $23 \Z[\sqrt{3}]=\fq \overline{\fq}$.
Since $\fq$ has residue field $\F_{23}$, we see that $\gamma^{11} \equiv \pm 1
\mod{\fq}$. Moreover, as $\tau(p) \equiv 2 \mod{23}$ we have
\[
2+\sqrt{3} \equiv \pm (2+\sqrt{3})^a \mod{\fq}.
\]
However, $2+\sqrt{3}$ has multiplicative order $11$ in $\Z[\sqrt{3}]/\fq=\F_{23}$.
As $0 \le a \le 10,$ we conclude that $a=1$. Thus
\[
\tau(p)+\sqrt{3}=(2+\sqrt{3}) (U+V\sqrt{3})^{11}, \qquad U,~V \in \Z.
\]
Comparing coefficients of $\sqrt{3}$ we obtain the Thue equation
\begin{multline*}
U^{11} + 22 U^{10} V + 165 U^9 V^2 + 990 U^8 V^3 + 2970 U^7 V^4 + 8316 U^6 V^5 +\\
        12474 U^5 V^6 + 17820 U^4 V^7 + 13365 U^3 V^8 + 8910 U^2 V^9 +
        2673 U V^{10} + 486 V^{11} \; =\; 1.
\end{multline*}
The \texttt{Magma} Thue equation solver (based
on algorithms in \cite{Smart}) gives that
the only solution is $(U,V)=(1,0)$.
Thus  $p=U^2-3V^2=1$ which is a contradiction.
\end{proof}

\noindent \textbf{Remark.}
The reader might be wondering if the case $11 \mid b$ can
also be tackled using Lemma~\ref{lem:sieve} instead of
appealing to Darmon and Merel. In that case, $N^\prime=32$,
and there is precisely one weight $2$ newform $f$ of level $32$.
This has rational eigenvalues and corresponds to the elliptic curve
\[
E \; : \; Y^2=X^3-X.
\]
Let $\ell \ne 2$ be a prime.
By inspection of the definition of $\cB_\ell$ in Lemma~\ref{lem:presieve},
we note that $(-1,0) \in \cB_\ell$ and that
$E_{-1,0,3}$ is the reduction modulo $\ell$
of the elliptic curve $E$. Thus $a_\ell(E_{-1,0,3})=a_\ell(E)=c_\ell$
where $c_\ell$ is the $\ell$-th coefficient of $f$. Thus
$(-1,0) \in \cC_{\ell,3}(f)$, and therefore $1 \in \cD_{\ell,3}(f)$.
Let $\kappa=1$. Going through the definitions in Lemma~\ref{lem:sieve},
it it easy to verify that $0 \in \cH_3(f)$ regardless of the choice
of $M$ and $\cL$. Hence we cannot use Lemma~\ref{lem:sieve}
to rule out the case $\kappa=1$ and $11 \mid b$.

\bigskip

There is a similar explanation for why we are unable to
use Lemma~\ref{lem:sieve} on its own to rule out the case
$q=3$, $\kappa=1$ and $11 \nmid b$. Here $N^\prime=96$.
There are two weight $2$ newforms of level $96$ and we
take $f$ to be the one corresponding to the elliptic curve
\[
E \; : \; Y^2=X^3+4X^2+3X.
\]
Let $\ell \nmid 6$ be  a prime.
We note that $(1,2) \in \cB_\ell$. Moreover, $E_{1,2,1}$ is
the reduction modulo $\ell$ of $E$. Hence $a_\ell(E_{1,2,1})=a_\ell(E)
=c_\ell$ which is as before the $\ell$-th coefficient of $f$.
We therefore have $(1,2) \in \cC_\ell(f)$ and so $3 \in \cD_{\ell,1}(f)$.
It follows, for $\kappa=1$, that $1 \in \cH_1(f)$
regardless of the choice of $M$ and $\cL$.

\begin{proof}[Proof of Proposition~\ref{prop:tpsmooth}]
Again we checked that equation \eqref{eqn:tpsmooth} has no
solutions with $p<200$ so we may suppose that $p>200$.
Moreover, by Lemmas \ref{lem:mod7}, \ref{lem:mod5} and \ref{lem:mod9},
we have $b_1=0$ or $b_1=1$, and $b_2=b_3=0$ in \eqref{eqn:tpsmooth}.
If $b_1=0$ then equation \eqref{eqn:tpsmooth} becomes
$\tau(p)^2-p^{11}=\pm 11^{b_4}$ which does not have any solutions
by Proposition~\ref{prop:pmqb}. Hence $b_1=1$. For convenience
we write $b$ for $b_4$, so equation \eqref{eqn:tpsmooth} becomes
\begin{equation}\label{eqn:BS2}
\tau(p)^2-p^{11}=\pm 3 \cdot 11^b.
\end{equation}
We apply Lemma~\ref{lem:sieve} with $q=11$
and $\kappa=\pm 3$. Here $N^\prime=96$ if $11 \mid b$
and $N^\prime=96 \times 11=1056$ if $11 \nmid b$.
For the newforms $f$ at both these levels and for
$\kappa=3$ and $\kappa=-3$, we computed $\cH_1(f)$ and $\cH_3(f)$.
We found that all these are empty with precisely one exception.
For that exception $\kappa=3$, and $f$ is the newform of level $96$
corresponding to the elliptic curve $E$ with Cremona label
\texttt{96a1}, when we find
\[
\cH_1(f)=\{0,\, 22,\, 44,\, 66,\, 88,\, 110,\, 132,\, 154,\, 176,\,
198,\, 220,\, 242,\, 264,\, 286,\, 308,\, 330,\,
    352,\, 374\},
\]
and so Lemma~\ref{lem:sieve} does not provide a contradiction.
However, we know that if $(p,b)$ is a solution to \eqref{eqn:BS2}
then $\overline{\rho}_{E_p,11} \sim \overline{\rho}_{f,\varpi} \sim
\overline{\rho}_{E,11}$. Suppose $b \ne 0$. Then the Frey-Hellegouarch curve
$E_p$ has conductor $96 \cdot 11$ and so multiplicative reduction
at $11$. The curve $E$ has conductor $96$ and hence good reduction
at $11$. Comparing the traces of Frobenius at $11$
in the two representations $\overline{\rho}_{E_p,11} \sim
\overline{\rho}_{E,11}$ (see \cite{KO}) we obtain
$\pm (11+1) \equiv a_{11}(E) \mod{11}$. However, $a_{11}(E)=4$
giving a contradiction. Thus $b=0$. Equation~\eqref{eqn:BS2}
now becomes equation~\eqref{eqn:pm3}, which we showed, in the proof of
Proposition~\ref{prop:pmqb}, to have no solutions.
This completes the proof.
\end{proof}

\section{The equation $\tau(p^4)= \kappa \cdot q^b$}\label{sec:taup4}

In this section, we establish the following two propositions.
\begin{prop}\label{prop:4pmqb}
Let $3 \le q <100$ be a prime.  The equation
\begin{equation}\label{eqn:4pmqb}
\tau(p^4) \; = \; \pm q^b,
\qquad
\text{$p$ prime,}
\qquad \text{$b \ge 0$}
\end{equation}
has no solutions.
\end{prop}
\begin{prop}\label{prop:4tpsmooth}
The equation
\begin{equation}\label{eqn:4tpsmooth}
\tau(p^4) \; = \; \pm 3^{b_1} 5^{b_2} 7^{b_3} 11^{b_4},
\qquad
\text{$p$ prime,}
\qquad \text{$b_1,~b_2,~b_3,~b_4 \ge 0$}
\end{equation}
has no solutions.
\end{prop}

We consider first the following general equation.
\begin{equation}\label{eqn:4quad}
\tau(p^4) \; = \; \kappa \cdot q^b,
\qquad
\text{$p \nmid 2\kappa q$ prime,}
\qquad \text{$b \ge 0$}.
\end{equation}
Here $\kappa$ is an odd integer, $q$ is an odd prime, and we assume for convenience
that
\[
q \nmid 5 \kappa, \qquad \text{$\ord_5(\kappa)=0$ or $1$}.
\]
Using the recursion \eqref{eqn:recursion} we find that
\[
\tau(p^4)=\tau(p)^4-3p^{11} \tau(p)^2 +p^{22}.
\]
which can be written as
\begin{equation}\label{eqn:rewrite}
4 \tau(p^4) = (2 \tau(p)^2-3 p^{11})^2 - 5 p^{22}.
\end{equation}
We may therefore rewrite \eqref{eqn:4quad} as
\[
5 (p^2)^{11}+ (4 \cdot \kappa \cdot q^b) \cdot 1^{11}\; =\; (2\tau(p)^2-3 p^{11})^2,
\]
which is an equation of signature $(11,11,2)$.
As before we follow the first author and Skinner \cite{BenS},
and associate to a solution of
\eqref{eqn:4quad} the Frey-Hellegouarch curve
\[
\begin{cases}
E_{p}\; : \;
Y^2=X(X^2+(3p^{11}-2 \tau(p)^2) X + \tau(p)^4-3p^{11} \tau(p)^2 +p^{22}) & \text{if $p \equiv 1 \mod{4}$},\\
E_{p} \; : \;
Y^2=X(X^2+(2 \tau(p)^2-3p^{11}) X +\tau(p)^4-3p^{11} \tau(p)^2 +p^{22}) & \text{if $p \equiv 3 \mod{4}$}.\\
\end{cases}
\]
Let
\begin{equation}\label{eqn:4N}
N=\begin{cases}
2^3 \cdot 5 \cdot \Rad(\kappa) \cdot q \cdot p & \text{if $b >0$, $\ord_5(\kappa)=0$}\\
2^3 \cdot 5 \cdot \Rad(\kappa) \cdot p & \text{if $b=0$, $\ord_5(\kappa)=0$}\\
2^3 \cdot 5^2 \cdot \Rad(\kappa/5) \cdot q \cdot p & \text{if $b >0$, $\ord_5(\kappa)=1$}\\
2^3 \cdot 5^2 \cdot \Rad(\kappa/5) \cdot p & \text{if $b=0$,
$\ord_5(\kappa)=1$},
\end{cases}
\end{equation}
and
\begin{equation}\label{eqn:4Nd}
N^\prime=\begin{cases}
2^3 \cdot 5 \cdot \Rad(\kappa) \cdot q  & \text{if $11 \nmid b$, $\ord_5(\kappa)=0$}\\
2^3 \cdot 5 \cdot \Rad(\kappa)  & \text{if $11 \mid b$, $\ord_5(\kappa)=0$}\\
2^3 \cdot 5^2 \cdot \Rad(\kappa/5) \cdot q  & \text{if $11 \nmid b$, $\ord_5(\kappa)=1$}\\
2^3 \cdot 5^2 \cdot \Rad(\kappa/5)  & \text{if $11 \mid b$,
$\ord_5(\kappa)=1$}.
\end{cases}
\end{equation}
The Frey curve $E_p$ has conductor $N$, and again it follows from
the recipes of the first author and Skinner \cite{BenS}
that there is a normalized newform $f$ as in \eqref{eqn:newform}
of weight $2$ and level $N^\prime$
and a prime $\varpi \mid 11$ in the integers of $K=\Q(c_1,c_2,\dotsc)$
so that \eqref{eqn:mod11} holds.

Throughout what follows, $\ell$ will be a prime satisfying
\begin{equation}\label{eqn:lassump4}
\ell \nmid 2 \cdot 5 \cdot 11 \cdot \kappa \cdot q \cdot p.
\end{equation}
As before \eqref{eqn:mod11norm} holds.

\begin{lem}\label{lem:presieve4}
Let $(p,b)$ be a solution to \eqref{eqn:4quad} and suppose $p \ne 3$, $23$.
Let $\ell$ be a prime satisfying \eqref{eqn:lassump4}. Let
\[
\cA_{\ell}=\{(s,t) \; : \; s,~t \in \F_\ell, \quad s \not \equiv 0 \mod{\ell},
\qquad
t^4-3s^{11} t^2 +s^{22} \not \equiv 0 \mod{\ell}
\},
\]
and let $\cB_{\ell}$ be as in \eqref{eqn:cB}.
For $(s,t) \in \cB_\ell$ let
\[
E_{s,t,1}/\F_\ell \; : \; Y^2=X(X^2+(3s^{11}-2t^2) X +t^4-3s^{11}t+s^{22}),
\]
\[
E_{s,t,3}/\F_\ell \; : \; Y^2=X(X^2+(2t^2-3s^{11})X +t^4-3s^{11}t+s^{22}).
\]
Let $f$ be a newform of weight $2$ and level $N^\prime$
so that \eqref{eqn:mod11} is satisfied, and $c_\ell$
be its $\ell$-th coefficient.
For $j=1$, $3$, let $\cC_{\ell,j}(f)$
be as in \eqref{eqn:cC}, and let
\[
\cD_{\ell,j}(f)=\{ t^4-3 s^{11} t^2+s^{22} \; : \; (s,t) \in \cC_{\ell,j}(f)\} \subseteq \F_\ell.
\]
If $p \equiv 1 \mod{4}$ then
$(\kappa \cdot q^b \mod \ell) \in \cD_{\ell,1}(f)$.
If $p \equiv 3 \mod{4}$ then
$(\kappa \cdot q^b \mod \ell) \in \cD_{\ell,3}(f)$.
\end{lem}
\begin{proof}
The proof is practically identical to that of Lemma~\ref{lem:presieve}.
\end{proof}

\begin{lem}\label{lem:sieve4}
Let $(p,b)$ be a solution to \eqref{eqn:4quad} with $p \ne 3$, $23$.
Let $M$ be a positive integer satisfying $22 \mid M$.
Let
\[
\cE=
\{ 0 \le \beta \le M-1 \; : \; \kappa \cdot q^\beta \equiv 1 \mod{4}\}.
\]
Let $f$ be a newform of weight $2$ and level $N^\prime$
so that \eqref{eqn:mod11} is satisfied.
Let
\[
\cF(f)=
\begin{cases}
\{ \beta \in \cE \; : \; 11 \nmid \beta\}
& \text{if $N^\prime=
2^3 \cdot 5 \cdot \Rad(\kappa) \cdot q$ or $2^3 \cdot 5^2 \cdot \Rad(\kappa/5) \cdot q$}\\
\{ \beta \in \cE \; : \; 11 \mid \beta\}
& \text{if $N^\prime=
2^3 \cdot 5 \cdot \Rad(\kappa)$ or $2^3 \cdot 5^2 \cdot \Rad(\kappa/5)$}\\
\end{cases}
\]
Let $\cL$ be a set of
primes satisfying
\begin{equation}\label{eqn:cL4}
\ell \nmid 2 \cdot 5 \cdot 11 \cdot \kappa q p, \qquad
O_{\ell}(q) \mid M.
\end{equation}
For $\ell \in \cL$ and $j=1$, $3$ 
\[
\cG_{\ell,j}(f) \; = \; \{ \beta \in \cF(f) \; : \;
(\kappa \cdot q^\beta \bmod \ell) \in \cD_{\ell,j}(f)
\}.
\]
Let
\[
\cH_j(f) \; = \; \bigcap_{\ell \in \cL} \cG_{\ell,j}(f).
\]
If $p \equiv 1 \mod{4}$ then
there is some $\beta \in \cH_1(f)$ such that
$b \equiv \beta \mod{M}$.
If $p \equiv 3 \mod{4}$ then
there is some $\beta \in \cH_3(f)$ such that
$b \equiv \beta \mod{M}$.
\end{lem}
\begin{proof}
This is almost identical to the proof of Lemma~\ref{lem:sieve}.
The main difference is that the sets $\cE$, $\cF(f)$ do not depend on the class
of $p$ modulo $4$, and we explain this now. Observe from \eqref{mod2}
that
\[
\kappa \cdot q^b = \tau(p^4) \equiv p^{44}+p^{33}+p^{22}+p^{11}+1 \equiv 3+2p \equiv 1 \pmod{4}
\]
regardless of the residue class of $p$ modulo $4$.
\end{proof}

\begin{lem}\label{lem:bne0}
The equations $\tau(p^4)=\pm 1$ and $\tau(p^4)=\pm 5$ have no solutions with $p$ prime.
\end{lem}
\begin{proof}
From the proof of Lemma~\ref{lem:sieve4}, we know that $\tau(p^4) \equiv 1
\pmod{4}$. Thus we need only consider the equations $\tau(p^4)=1$ and
$\tau(p^4)=5$.
Suppose $\tau(p^4)=1$ and write $z=2 \tau(p)^2-3p^{11}$. From \eqref{eqn:rewrite} we have
\[
z^2-5 p^{22}= 4.
\]
Write $\varepsilon=(1+\sqrt{5})/2$.
Then $(\lvert z\rvert +p^{11} \sqrt{5})/2$ is
a positive unit in $\Z[\varepsilon]$ with norm $+1$.
Hence
\[
\frac{z+p^{11} \sqrt{5}}{2} \; = \; \varepsilon^{2n}, \qquad
\varepsilon=(1+\sqrt{5})/2.
\]
for some $n \in \Z$. Thus
\[
p^{11}= \frac{\varepsilon^{2n}-\overline{\varepsilon}^{2n}}{\sqrt{5}}
=F_{2n}
\]
where $F_n$ denotes
the $n$-th Fibonacci number.
By \cite{BMS}
the only perfect powers in the Fibonacci sequence are
$0$, $1$, $8$ and $144$, giving a contradiction. Alternatively,
$F_{2n}=F_n L_n$ where $L_n$ is the $n$-th Lucas number.
From the identity $L_n^2-5 F_n^2=4 \cdot (-1)^n$ we
see that $\gcd(F_n,L_n)=1$ or $2$. Thus $F_n=1$ or $L_n=1$
quickly leading to a contradiction.

Next we suppose  that  $\tau(p^4)=5$ and 
write $z=5w$. Hence
\[
5w^2-p^{22}=4
\]
and it follows that there is an integer $n$ such that
\[
p^{11}=\varepsilon^n+\overline{\varepsilon}^n=L_n,
\]
where $L_n$ denotes the $n$-th Lucas number.
By \cite{BMS},
the only perfect powers in the Lucas sequence are $1$ and $4$,
again giving a contradiction.
\end{proof}

\begin{proof}[Proof of Proposition~\ref{prop:4pmqb}]
We checked that \eqref{eqn:4pmqb} has no solutions
with $p<200$ for primes $3 \le q <100$. We shall henceforth suppose
that $p>200$. In particular, $p \ne q$.
Moreover, any solution to \eqref{eqn:4pmqb}
is a solution to \eqref{eqn:4quad} with $\kappa=1$ or $-1$.

We consider $q=5$ first. By \eqref{eqn:rewrite},
$\ord_5(\tau(p^4))=0$ or $1$. Thus we reduce to the equations
$\tau(p^4)=\pm 1$ and $\tau(p^4)=\pm 5$. These do not have solutions
by Lemma~\ref{lem:bne0} and hence we may assume that $q \ne 5$. From Lemma~\ref{lem:bne0} again
we have $b>0$. By \eqref{eqn:rewrite},
$5$ is a quadratic residue modulo $q$. The possible
values of $q$ are
\begin{equation}\label{eqn:q}
11,\; 19,\; 29,\; 31,\; 41,\; 59,\; 61,\; 71,\; 79,\; 89.
\end{equation}
For each of these values we take
\begin{equation}\label{eqn:McL}
M=396=2^2 \cdot 3^2 \cdot 11, \qquad
\cL=\{\text{$3 \le \ell <200$ prime, $\ell \ne 5$, $11$, $q$} \; : \;  O_\ell(q) \mid M\}.
\end{equation}
Observe that since $p>200$ that $\ell \ne p$, and thus
satisfies \eqref{eqn:cL4}.

We consider first the case $11 \nmid b$.
Thus, in \eqref{eqn:4Nd}, the level $N^\prime=2^3 \cdot 5 \cdot q$.
We computed for each newform $f$ of level $N^\prime$
the sets $\cH_1(f)$ and $\cH_3(f)$, both for $\kappa=1$, $\kappa=-1$.
We found all of these
to be empty. By Lemma~\ref{lem:sieve}, we conclude that
\eqref{eqn:4pmqb} has no solutions with $11 \nmid b$.

Next we consider $11 \mid b$. Thus $N^\prime=2^3 \cdot 5$.
There is a unique newform $f$ of level $N^\prime$ which corresponds
to the elliptic curve $E$ with Cremona label \texttt{40a1}.
Thus, from \eqref{eqn:mod11}
we obtain $\overline{\rho}_{E_p} \sim \overline{\rho}_{E}$.
Note, by \eqref{eqn:4N}
that $E_p$ has multiplicative reduction at $q$. However, $E$ has good reduction
	at $q$. Thus, by \cite{KO}, 
	we have $\pm (q+1) \equiv a_q(E) \pmod{11}$. We checked that this
does not hold for all the values of $q$ in \eqref{eqn:q}.
This completes the proof.
\end{proof}

\begin{proof}[Proof of Proposition~\ref{prop:4tpsmooth}]
Again we checked that \eqref{eqn:4tpsmooth} has no solutions with $p<200$, whence
we may suppose $p>200$. Moreover, as $5$ is a quadratic non-residue modulo
$3$ and $7$, we see from \eqref{eqn:rewrite} that $b_1=b_3=0$ in
\eqref{eqn:4tpsmooth}. Also $5^2 \nmid \tau(p^4)$ from \eqref{eqn:rewrite},
so $b_2=0$ or $1$.
But from Proposition~\ref{prop:4pmqb} we have $b_2 \ne 0$, and so $b_2=1$.
We have thus reduced to consideration of the equation
\[
\tau(p^4)=\pm 5 \cdot 11^b,
\]
whereby we have $\kappa=\pm 5$ and $q=11$.
Observe that $b>0$ by Lemma~\ref{lem:bne0}. Suppose $11 \nmid b$.
Thus $N^\prime=8 \cdot 25 \cdot 11=2200$.
We take $M$ and $\cL$ as in \eqref{eqn:McL}.
There are $25$ conjugacy
classes of newforms $f$ of weight $2$ and level $2200$. For each, we
found $\cH_1(f)$ and $\cH_3(f)$ to be empty, both for $\kappa=5$ and $\kappa=-5$.
By Lemma~\ref{lem:sieve4}, there are no solutions with $11 \nmid b$.
Thus $11 \mid b$, and so $N^\prime=2^3 \cdot 25=200$. There are
five weight $2$ newforms of level $200$. We computed $\cH_1(f)$
and $\cH_3(f)$ for these, both for $\kappa=5$ and $\kappa=-5$.
The only non-empty one we found was $\cH_3(f)$ for $\kappa=5$
where $f$ is the rational newform corresponding to the elliptic curve $E$
with Cremona label \texttt{200b1}.
Then $\overline{\rho}_{E_p,11} \sim \overline{\rho}_{E,11}$.
Here $E_p$ has multiplicative
reduction at $11$, though $E$ has good reduction at $11$.
As before,
$\pm (11+1) \equiv a_{11}(E) \mod{11}$.
However, $a_{11}(E)=-4$, giving a contradiction and completing the proof.
\end{proof}

\section{On the largest prime divisor of $\tau(p^3)$}\label{sec:taup3}
\begin{prop}\label{prop:tp3smooth}
Let $p$ be a prime for which $\tau(p) \ne 0$. Then
$P(\tau(p^3)) \ge 13$, unless $p=2$, in which case we have
$\tau(8)=2^9 \cdot 3 \cdot 5 \cdot 11$.
\end{prop}
We consider
\begin{equation}\label{eqn:tp3smooth}
P(\tau(p^3)) \le 11.
\end{equation}
We checked that the only $p< 200$ satisfying \eqref{eqn:tp3smooth}
is $p=2$.
We shall therefore suppose $p>200$.
Recall that $\tau(p^3)=\tau(p) \cdot (\tau(p)^2-2 p^{11})$.
From \eqref{mod3} and \eqref{mod7}, we easily see that $3$ and $7$
do not divide $\tau(p)^2-2p^{11}$. Moreover, we recall that $\tau(p)$
is even, so $\ord_{2}(\tau(p)^2-2p^{11})=1$. Thus
\begin{equation}\label{eqn:factored}
\tau(p)^2-2 p^{11}= \pm 2 \cdot 5^a \cdot 11^b \; \; \mbox{ and } \; \;  \tau(p)=\pm 2^r \cdot 3^s
\cdot 5^t \cdot 7^u \cdot 11^v.
\end{equation}
As before, we associate to this a Frey-Hellegouarch curve
\[
E_p \; : \; Y^2=X(X^2+2 \tau(p) X + 2 p^{11}).
\]
By the recipes of the first author and Skinner, 
the conductor of $E_p$ is one of
\[
N=2^8 \cdot p, \qquad 2^8\cdot 5 \cdot p, \qquad 2^8 \cdot 11 \cdot p, \qquad 2^8 \cdot 5 \cdot
11 \cdot p,
\]
and \eqref{eqn:mod11}
holds for some weight $2$ newform $f$ whose level $N^\prime$
is one of the following
\begin{equation}\label{eqn:Ndhere}
N^\prime=2^8, \qquad 2^8\cdot 5, \qquad 2^8 \cdot 11, \qquad 2^8 \cdot 5 \cdot
11.
\end{equation}
There are a total of $123$ conjugacy classes of newforms $f$ at these levels.
Let $f$ be any of these such that \eqref{eqn:mod11} holds. Let $\ell \ne 2$,
$5$, $11$, $p$ be a prime. Then $11 \mid \Norm(a_\ell(E_p)-c_\ell(f))$ where $c_\ell(f)$
is the $\ell$-th coefficient of $f$.
\begin{lem}
Let $\ell \ne 2$, $5$, $11$ be a prime $<200$.
Let $p$ be an odd prime with $\tau(p) \ne 0$ and $P(\tau(p^3)) \leq 11$.
Let $f$ be a newform of weight $2$ and one of the levels $N^\prime$
in \eqref{eqn:Ndhere} so that \eqref{eqn:mod11} is satisfied.
Write
\[
\cA_{\ell}=
\begin{cases}
\{(s,t) \; : \; s,~t \in \F_\ell, \quad s (t^2-2s^{11})\not \equiv 0
\mod{\ell}\}, \qquad \ell=3,~7\\
\{(s,t) \; : \; s,~t \in \F_\ell, \quad s t (t^2-2s^{11})\not \equiv 0
\mod{\ell}\}, \qquad \ell \ge 13.\\
\end{cases}
\]
Let $\cB_\ell$ be as in \eqref{eqn:cB}. Let
\[
E_{s,t}/\F_{\ell} \; : \; Y^2=X(X^2+2tX+s^{11}),
\]
and
\[
\cC_{\ell}(f)=\{(s,t) \in \cB_\ell \; : \; \Norm(a_\ell(E_{s,t})-c_\ell(f))
\equiv 0 \pmod{11}\}.
\]
Then there is some $(s,t) \in \cC_\ell(f)$ so that
$(p,\tau(p)) \equiv (s,t) \pmod{\ell}$.
\end{lem}
\begin{proof}
This is is similar to the proof of Lemma~\ref{lem:presieve}.
\end{proof}
\begin{proof}[Proof of Proposition~\ref{prop:tp3smooth}]
For each of the $123$ conjugacy classes of newforms $f$ we computed
$\cC_{\ell}(f)$ for $\ell=3$, $7$, $13$ and $23$. We found that
at least one of these four empty, except for the three rational
newforms which correspond to the elliptic curves \texttt{256a1},
\texttt{256b1} and \texttt{256c1}. All three elliptic curves have
CM, respectively by $\Q(\sqrt{-2})$, $\Q(\sqrt{-1})$, $\Q(\sqrt{-1})$.
Note that $11$ splits in $\Q(\sqrt{-2})$ and is inert in $\Q(\sqrt{-1})$.
Hence the image of $\rho_{E_p,11} \sim \rho_{f,11}$ belongs
to the normalizer of split Cartan subgroup in the first case,
and the normalizer of a non-split Cartan subgroup in the second
and third case. Thanks to the work of Momose \cite{Momose},
and Darmon and Merel \cite[Theorem 8.1]{DM}, the $j(E_p) \in \Z[1/11]$.
However, $E_p$ has multiplicative reduction at $p$ giving a contradiction.
\end{proof}

\section{Proof of Theorem \ref{thm2a}}\label{sec:thm2a}

\begin{lem}\label{lem:ge13}
Let $p \le 11$ be a prime. Suppose
$P(\tau(p^{m-1})) \le 11$ with $m \ge 3$.
Then $p=2$ and $m=4$.
\end{lem}
\begin{proof}
First let $p=2$.
Let $m \ge 3$ be such that $P(\tau(2^{m-1})) \le 11$.
Note that $\tau(2)=-2^3 \times 3$.
Let $\{u_n\}$ be the Lucas sequence defined in Lemma~\ref{lem:associated},
with characteristic polynomial $X^2-3X+2^5$.
Then $P(u_m) \le 11$. Moreover, by part (i) of Theorem~\ref{thm:Carmichael},
we have $2 \nmid u_n$ for all $n \ge 1$.
We note that
\[
u_2=-3, \quad u_3=-23, \quad
u_4=3 \times 5 \times 11, \quad u_5=241,
\]
\[ u_6=-3^2 \times 23 \times 29,
\quad
u_7=7 \times 1471, \quad u_8=3 \times 5 \times 11 \times 977.
\]
By the Primitive Divisor Theorem (Theorem~\ref{thm:BHV}),
every term $u_n$ with $n \ge 9$ is divisible
by some prime $\ell \ge 13$. Thus the only terms with $P(u_n) \le 11$
are $u_2$ and $u_4$. Since $m \ge 3$ we have $m=4$.

By a similar strategy we checked that $P(\tau(p^{m-1})) \ge 13$
for $3 \le p \le 11$ and $m\ge 3$.
\end{proof}

\begin{lem}\label{lem:taupm}
Let $p$ be a prime. Let $m \ge 3$ be an integer such
that $\tau(p^{m-1}) \ne 0$ and
\begin{equation}\label{taupm}
P(\tau (p^{m-1})) \leq 11.
\end{equation}
Then $p=2$ and $m=4$.
\end{lem}
\begin{proof}
By Lemma~\ref{lem:ge13},  we may suppose $p \ge 13$.
If $\tau(p)=0$, by Lemma~\ref{taup=0}, we have $\tau(p^{m-1})=0$
or a power of $p$ contradicting the hypotheses of the lemma.
We may therefore suppose $\tau(p) \ne 0$.
Fix $p$ and let $m$ be the least value $ \ge 3$ such that \eqref{taupm}
is satisfied.
By Propositions~\ref{prop:tpsmooth}, \ref{prop:4tpsmooth}
and \ref{prop:tp3smooth}, we know that $m \ge 6$.

Suppose first that $p \mid \tau(p)$. By induction from
\eqref{eqn:recursion} we have $p \mid \tau(p^n)$ for all $n \ge 1$.
Hence $p \mid \tau(p^{m-1})$ and so $p \le 11$ giving a contradiction.
Thus $p \nmid \tau(p)$. Let $u_n=\tau(p^{n-1})$ for $n \ge 1$.
Then $\{u_n\}$ is a Lucas sequence by Lemma~\ref{lem:associated}.
Now $u_k \mid u_n$ if $k \mid n$. As $m \ge 6$, it is divisible
either by $4$ or an odd prime. However $u_4=\tau(p^3)$,
and so $P(u_4) \ge 13$ by Proposition~\ref{prop:tp3smooth}.
Hence $m$ is divisible by an odd prime, and from the minimality
of $m$ it follows that $m \ge 7$ is a prime. 
By the Primitive Divisor Theorem, $u_m=\tau(p^{m-1})$ has a prime divisor $q$
that does not divide
$u_1 u_2 \cdots u_{m-1}$ nor $D=(\alpha-\beta)^2$
(where $\alpha$, $\beta$ are as in Lemma~\ref{lem:associated}).
Here $q=2$, $3$, $5$, $7$ or $11$. 
But $m_q$, the rank of apparition of $q$, divides $m$ by
Theorem~\ref{thm:Carmichael}, and so $m_q=m$. However
$m_q \mid (q-1)(q+1)$, again from Theorem~\ref{thm:Carmichael}.
But $(q-1)(q+1)$ is not divisible
by a prime $\ge 7$ for $q=2$, $3$, $5$, $7$ or $11$.
This contradiction completes the proof.
\end{proof}

\begin{proof}[Proof of Theorem~\ref{thm2a}]
Suppose $n$
is a powerful number such that $\tau(n) \ne 0$ and $P(\tau(n)) \le 11$.
Let $p$ be a prime divisor of $n$.
Thus $p^{m-1} \mid\mid n$, where, as $n$ is powerful,
$m \ge 3$. Now, as $\tau$ is multiplicative,
$\tau(p^{m-1}) \ne 0$ and
$\tau(p^{m-1}) \mid \tau(n)$.
In particular, $P(\tau(p^{m-1})) \le 11$.
By Lemma~\ref{lem:taupm} we have $p=2$ and $m=4$.
Thus $n=8$ as required.
\end{proof}

\section{Proof of Theorem~\ref{thm3}}\label{sec:thm3}

\begin{proof}[Proof of Theorem~\ref{thm3}]
Suppose $\tau(n)=\pm q^\alpha$ where $3 \le q <100$
is a prime and $n \ge 2$. Then $\tau(n)$
is odd, and so $n$ must be an odd square. Thus there is
a prime $p$ and an integer $m \ge 3$
such that $p^{m-1} \mid\mid n$. Hence $\tau(p^{m-1})=\pm q^a$
for some $a \ge 0$. The following lemma completes the proof.
\end{proof}

\begin{lem}
Let $3 \le q <100$ be a prime and $a$ be a nonnegative integer. Let $p$ be a prime and $m \ge 3$
an odd integer. Then $\tau(p^{m-1}) \ne \pm q^a$.
\end{lem}
\begin{proof}
	We argue by contradiction. 
Suppose $m \ge 3$ is the smallest
odd integer such that
\begin{equation}
\tau(p^{m-1})=\pm q^a.
\end{equation}
By Propositions~\ref{prop:pmqb} and ~\ref{prop:4pmqb},
we have $m \ge 7$.
We treat first the case $\tau(p)=0$.
By Lemma~\ref{taup=0},
we see that $\pm q^a=\tau(p^{m-1})$ is either $0$ or a power of $p$.
Thus $p=q<100$, which gives a contradiction
since any $p$ for which $\tau(p)=0$ satisfies
 \eqref{eqn:Derickx}.
Thus $\tau(p) \ne 0$.

Let $\{u_n\}$ be the 
Lucas sequence given in Lemma~\ref{lem:associated}.
It follows from that lemma that $u_n \mid \tau(p^{n-1})$
and $p \nmid u_n$ for all $n \ge 1$.
If $p =q$, then $u_{m}= \pm 1$
which contradicts the Primitive Divisor
Theorem (Theorem~\ref{thm:BHV}), as $m \ge 7$. We conclude that
$p \ne q$.

Next we consider the case $p \mid \tau(p)$.
Then $p \mid \tau(p^n)$ for all $n \ge 1$
by \eqref{eqn:recursion}, and so $p=q$, giving a contradiction.
Thus $p \nmid \tau(p)$.
It follows that $u_n=\tau(p^{n-1})$ for all $n \ge 1$.
Recall that if $k \mid n$ then $u_k \mid u_n$.
By the minimality of $m$ we see that $m \ge 7$ is a prime.
We invoke the Primitive Divisor Theorem again to conclude
that $q \nmid (\alpha-\beta)^2 u_1 u_2 \cdots u_{m-1}$
(in the notation of Lemma~\ref{lem:associated}).
From Theorem~\ref{thm:Carmichael}, $m=m_q \mid (q-1)(q+1)$.
The possible pairs of primes
$(q,m)$ with $3 \le q < 100$ and $m \mid (q-1)(q+1)$ are
$$
\begin{array}{l}
(13, 7), \;
(23, 11), \;
(29, 7), \;
(37, 19), \;
(41, 7),\;
(43, 7),\;
(43, 11), \;
(47, 23), \;
(53, 13),\;
(59, 29),\\
(61, 31), \;
(67, 11),\;
(67, 17),\;
(71, 7), \;
(73, 37),\;
(79, 13),\;
(83, 7),\;
(83, 41),\;
(89, 11),\;
(97, 7).
\end{array}
$$
Fixing any of these pairs $(q,m)$, it remains to solve
$\tau(p^{m-1})=\pm q^a$.
By Lemma~\ref{lem:tauFm}, and the fact that $m$ is prime, 
we see that $(X,Y,a)=(p^{11},\tau(p),a)$
is a solution to the Thue--Mahler equation
\[
	\Psi_m(X,Y)=\pm q^{a}.
\]
We solved these Thue--Mahler equations
using the \texttt{Magma} implementation
of the Thue--Mahler solver described in \cite{GKMS}.
None of the solutions are of the form $(p^{11},\tau(p),a)$.
This completes the proof of Theorem~\ref{thm3}. We illustrate this by providing
more details for  the case $q=83$. Here $m$ is a prime $\ge 7$ dividing $83^2-1=2^3 \times 3
\times 7 \times 41$, and thus the possible pairs $(q,m)$ are $(83,7)$ and
$(83,41)$. For the first pair, the Thue--Mahler equation is
\[
-X^3 + 6 X^2 Y - 5 X Y^2 + Y^3=83^a,
\]
and the solutions are
\[
(X,Y,a)= ( 5, 1, 0 ),\;
    ( -9, -14, 0 ),\;
    ( 2, 3, 0 ),\;
    ( -7, -1, 1 ),\;
    ( 5, 2, 1 ),\;
    ( 0, 1, 0 ),\;
    ( -1, -2, 0 ),
    \]
    \[
    ( -17, -38, 2 ),\;
    ( -8, -13, 1 ),\;
    ( 13, 20, 1 ),\;
    ( 1, 1, 0 ),\;
    ( 4, 13, 0 ),\;
    ( -6, -19, 1 ),\;
    ( -1, 0, 0 ),\;
    ( 21, 25, 2 ),
    \]
    \[
    ( 3, 11, 1 ),\;
    ( -4, 13, 2 ),\;
    ( -1, -3, 0 ),\;
    ( -5, -2, 1 ),\;
    ( 0, -1, 0 ),\;
    ( 17, 38, 2 ),\;
    ( 6, 19, 1 ),\;
    ( 7, 1, 1 ),
    \]
    \[
    ( 1, 0, 0 ),\;
    ( -4, -13, 0 ),\;
    ( 4, -13, 2 ),\;
    ( 9, 14, 0 ),\;
    ( -3, -11, 1 ),\;
    ( 1, 3, 0 ),\;
    ( -1, -1, 0 ),
    \]
    \[
    ( -13, -20, 1 ),\;
    ( -5, -1, 0 ),\;
    ( -21, -25, 2 ),\;
    ( 8, 13, 1 ),\;
    ( 1, 2, 0 ),\;
    ( -2, -3, 0 ).
\]
For the pair
$(q,m)=(83,41)$ the Thue--Mahler equation
is
\begin{multline*}
X^{20} - 210 X^{19} Y + 7315 X^{18} Y^2 - 100947 X^{17} Y^3 + 735471 X^{16} Y^4 -
    3268760 X^{15} Y^5\\ + 9657700 X^{14} Y^6 - 20058300 X^{13} Y^7 + 30421755
    X^{12} Y^8
    - 34597290 X^{11} Y^9 + 30045015 X^{10} Y^{10}\\ - 20160075 X^9 Y^{11} +
    10518300 X^8 Y^{12} - 4272048 X^7 Y^{13} + 1344904 X^6 Y^{14} - 324632 X^5
    Y^{15} \\
    +
    58905 X^4 Y^{16} - 7770 X^3 Y^{17} + 703 X^2 Y^{18} - 39 X Y^{19} + Y^{20}
= \pm 83^{a},
\end{multline*}
and the solutions are
\[
    ( -1, -3, 0 ),\;
    ( -1, -2, 0 ),\;
    ( 1, 2, 0 ),\;
    ( 1, 0, 0 ),\;
    ( -1, 0, 0 ),
\]
\[
    ( 1, 3, 0 ),\;
    ( 0, 1, 0 ),\;
    ( 0, -1, 0 ),\;
    ( 1, 1, 0 ),\;
    ( -1, -1, 0 ).
\]
\end{proof}

\noindent \textbf{Remark.}
The aforementioned Thue--Mahler solver requires
knowledge of the class group and unit group of
the number field defined by the equation $\Psi_{m}(1,Y)=0$;
this number field has degree $\phi(m)/2=(m-1)/2$. Ordinarily, if the
degree is too large, this might not be practical,
or might require assuming the Generalized Riemann Hypothesis. However,
from Lemma~\ref{lem:abelian}, this number field is 
$\Q(\zeta_m)^+$.
For the values of $m$ under consideration (and in fact
for all prime $m \le 67$), the class number $h_m^+$ of
$\Q(\zeta_m)^+$ is known to be $1$; see for example
\cite[Theorem 1]{Linden}. Moreover, if we denote
the unit group of $\Q(\zeta_m)^+$ by $E_m^+$ and the subgroup
of cyclotomic units by $C_m^+$ then $[E_m^+:C_m^+]=h_m^+$; see
\cite[Theorem 8.2]{Wa}. Hence in all our cases, $E_m^+=C_m^+$, and
is generated \cite[Lemma 8.1]{Wa} by $-1$ and $(1-\zeta_m^a)/(1-\zeta_m)$
with $1<a<m/2$. Thus for all values of $m$ under consideration
we know the class group and unit group.

\section{Concluding remarks}

As noted in the introduction, it would likely be extremely challenging computationally to extend, for example, Corollary \ref{thm2} to explicitly find all $n$ with $\tau (n)$ odd and, say,
$$
P(\tau (n)) \leq 17.
$$
The bottleneck in our approach is related to the difficulty involved in classifying the primes $p$ for which $P(\tau (p^2))$ and $P(\tau (p^4))$ are ``small''. 
For larger exponents $m$, finding the $p$ with $P(\tau (p^m))$ bounded appears to be somewhat more tractable.
By way of example,
we may show, by direct application of the Thue-Mahler solver described in  \cite{GKMS}, the following result.
\begin{prop}
The equation
\begin{equation}\label{eqn:6psmooth}
\tau(p^6) \; = \; \pm 3^{b_1} 5^{b_2} 7^{b_3} 11^{b_4} 13^{b_5} 17^{b_6} 19^{b_7} 23^{b_8 } 29^{b_9} 31^{b_{10}} 37^{b_{11}},
\qquad
\text{$p$ prime,}
\qquad \text{$b_i \in \mathbb{Z}$}
\end{equation}
has no solutions.
\end{prop}

This amounts to solving the Thue-Mahler equation
\[
-X^3 + 6 X^2 Y - 5 X Y^2 + Y^3= \pm 3^{b_1} 5^{b_2} 7^{b_3} 11^{b_4} 13^{b_5} 17^{b_6} 19^{b_7} 23^{b_8 } 29^{b_9} 31^{b_{10}} 37^{b_{11}}
\]
and checking to see if any solutions have $X=p^{11}$ for some prime $p$.
Appealing to  \cite{GKMS}, we find that all solutions in coprime integers $X$ and $Y$ have either $\max \{ |X|,|Y| \} < 100$, or satisfy
$$
\begin{array}{r} 
\pm (X,Y) \in  \left\{ (31,105),  (33,107),  (33,109), (41,124), (67,219), (74,115),  (74,117), (76,119),  \right. \\
 \left (83,125), (152,237), (207,-152), (251, 815), (313,62), (359,925), (564,877),  (566,773) \right\}. \\
\end{array}
$$


\begin{thebibliography}{}


\bibitem{BaCrOn}
J. S. Balakrishnan, W. Craig and K. Ono,
\newblock Variations of Lehmer's conjecture for Ramanujan's tau-function,
\newblock {\em J. Number Theory}, to appear.

\bibitem{BaCrOnTs}
J. S. Balakrishnan, W. Craig, K. Ono and W.-L. Tsai,
\newblock Variants of Lehmer's speculation for newforms,
\newblock submitted for publication.

\bibitem{BenHyp}
M. A. Bennett,
\newblock Rational approximation to algebraic numbers of small height : the Diophantine equation $|ax^n-by^n|=1$,
\newblock {\em J.\ Reine Angew.\ Math.} 535 (2001), 1--49.

\bibitem{BenS}
M.\ A.\ Bennett and C.\ Skinner,
\newblock Ternary Diophantine equations via Galois representations and modular forms,
\newblock {\em Canad.\ J.\ Math.} 56 (2004), no.\ 1, 23--54.

\bibitem{Barr}
C. F. Barros,
\newblock On the Lebesgue-Nagell equation and related subjects,
\newblock PhD thesis, University of Warwick, 2010.

\bibitem{BDMS}
M. A. Bennett, S. Dahmen, M. Mignotte and S. Siksek,
\newblock Shifted powers in binary recurrence sequences,
\newblock {\em Math. Proc. Cambridge Philos. Soc.}  158 (2015), 305--329.

\bibitem{BHV}
Y. Bilu, G. Hanrot and P. Voutier,
\newblock Existence of primitive divisors of Luca and Lehmer numbers,
\newblock {\em J.\ Reine Angew.\ Math.} 539 (2001), 75--122.

\bibitem{magma}
W.\ Bosma, J.\ Cannon and C.\ Playoust,
\newblock The Magma Algebra System I: The User Language,
\newblock {\em J.\ Symb.\ Comp.} {\bf 24} (1997), 235--265. (See also \url{http://magma.maths.usyd.edu.au/magma/})

\bibitem{Bug-BLMS}
Y. Bugeaud,
On the greatest prime factor of $ax^m-by^n$ II,
		\emph{Bull. London Math. Soc.}  32 (2000), no. 6, 673–678.

\bibitem{BuGy}
Y. Bugeaud and K. Gy\H{o}ry,
\newblock Bounds for the solutions of Thue-Mahler equations and norm form equations,
\newblock {\em Acta Arith.} 74 (1996), 273--292.






\bibitem{BMS}
Y. Bugeaud, M. Mignotte and S. Siksek,
\newblock Classical and modular approaches to exponential Diophantine equations I. Fibonacci and Lucas perfect powers,
\newblock {\em Ann. of Math.} 163 (2006), 969--1018.



\bibitem{Car}
R. D. Carmichael,
\newblock On the numerical factors of the arithmetic forms $\alpha^n \pm \beta^n$,
\newblock {\em Ann. Math.} 15 (1913), 30--70.


\bibitem{DM}
H.\ Darmon and L.\ Merel,
\newblock Winding quotients and some variants of Fermat's Last Theorem,
\newblock {\em J. Reine Angew.\ Math.} 490 (1997), 81--100.


\bibitem{Deligne}
P. Deligne,
\newblock La conjecture de Weil I,
\newblock {\em Inst. Hautes \'Etudes Sci. Publ. Math.} 43 (1974), 273--307.


\bibitem{DeJa}
S. Dembner and V. Jain,
\newblock Hyperelliptic curves and newform coefficients,
\newblock submitted for publication.

\bibitem{DeHoZe}
M. Derickx, M. Van Hoeij and Jinxiang Zeng,
\newblock Computing Galois representations and equations for modular curves $X_H(\ell)$,
\newblock arXiv:1312.6819v2.







\bibitem{GKMS}
A. Gherga, R. von K\"anel, B. Matschke and S. Siksek,
\newblock Efficient resolution of Thue--Mahler equations,
to appear.


\bibitem{HaMa}
M. Hanada and R. Madhukara,
\newblock Fourier coefficients of level $1$ Hecke eigenforms,
\newblock submitted for publication.


\bibitem{KO}
A. Kraus and J. Oesterl\'e,
\newblock Sur une question de B.\ Mazur,
\newblock {\em Math. Ann.} 293 (2002), 259--275.


\bibitem{Leh}
D. H. Lehmer,
\newblock The vanishing of Ramaujan's function $\tau (n)$,
\newblock {\em Duke Math. J.} 14 (1947), 429--433.

\bibitem{Lenstra}
H. W. Lenstra,
\newblock Algorithms in algebraic number theory,
\emph{Bull. Amer. Math. Soc.} \textbf{26} (1992), 211--244.

\bibitem{Linden}
F.\ J.\ van der Linden,
\newblock Class number computations of real abelian number fields,
\newblock {\em Mathematics of Computation} \textbf{39} (1982), 693--707.

\bibitem{LuMaSt}
F. Luca, S. Mabaso and P. Stanica,
\newblock On the prime factors of the iterates of the Ramanujan $\tau$-function,
\newblock {\em Proc. Edinburgh Math. Soc.}, to appear.


\bibitem{LuSh}
F. Luca and I. E. Shparlinski,
\newblock Arithmetic properties of the Ramanujan function,
\newblock {\em Proc. Indian Acad. Sci. (Math. Sci.)} 116 (2006), 1--8.


\bibitem{LyRo0}
N. Lygeros and O. Rozier,
\newblock  A new solution to the equation $\tau (p) \equiv 0 \mod{p}$,
\newblock {\em J. Integer Sequences} 13 (2010),  Article 10.7.4.

\bibitem{LyRo}
N. Lygeros and O. Rozier,
\newblock  Odd prime values of the Ramanujan tau function,
\newblock {\em Ramanujan Journal} 32 (2013),  269--280.

\bibitem{Momose}
F.\ Momose,
\newblock Isogenies of prime degree over number fields,
\newblock {\em Compositio Math.} 97 (1995), 329--348.

\bibitem{Mor}
L. J. Mordell,
\newblock On Mr. Ramanujan's empirical expansions of modular functions,
\newblock {\em Proc. Cambridge Phil. Soc.} 19 (1917), 117--124.

\bibitem{MuMu}
M. Ram Murty and V. Kumar Murty,
\newblock Odd values of Fourier coefficients of certain modular forms,
\newblock {\em Int. J. Number Theory} 3 (2007), 455-470.

\bibitem{MuMuSh}
M. Ram Murty, V. Kumar Murty and T. N. Shorey,
\newblock Odd values of the Ramanujan $\tau$-function,
\newblock {\em Bulletin Soc. Math. France} 115 (1987), 391--395.




\bibitem{OnTa}
K. Ono and Y. Taguchi,
\newblock $2$-adic properties of certain modular forms and their applications to arithmetic functions,
\newblock {\em Int. J. Number Theory} 1 (2005), 75--101.


\bibitem{Ram}
S. Ramanujan,
\newblock On certain arithmetical functions,
\newblock {\em Trans. Camb. Philos. Soc.} 22 (1916), 159--184.

\bibitem{Ser}
J.-P.\ Serre,
\newblock Divisibilit\'e de certaines fonction arithmetiques,
\newblock {\em L'Ens. Math.} 22 (176), 227--260.



\bibitem{Siksek}
S. Siksek,
\newblock The Modular Approach to Diophantine Equations,
in \emph{Explicit Methods in Number Theory: Rational Points and Diophantine
Equations}, ed. Belabas et al.,
Panoramas et synth\`{e}ses \textbf{36}, 2012.


\bibitem{Smart}
N.\ P.\ Smart,
\newblock The Algorithmic Resolution of Diophantine Equations,
\newblock Cambridge University Press, 1998.

\bibitem{St13}
C. L. Stewart,
\newblock On divisors of Lucas and Lehmer numbers,
\newblock {\em Acta Math.} 211 (2013), 291--314.


\bibitem{SwD}
H.\ P.\ F.\ Swinnerton-Dyer, 
\newblock On $\ell$-adic representations and congruences for coefficients of modular forms, 
page 1--55 of  W.\ Kuyk and J.-P.\ Serre (eds.), Modular functions of one variable, III, Lecture Notes in Mathematics 350, 1973. 

\bibitem{Wa}
L.\ C.\ Washington,
Introduction to Cyclotomic Fields, Springer Verlag,
1982.
\end{thebibliography}
\end{document}